\documentclass[12pt,thmsb,titlepage]{article}
\usepackage{amsmath}
\usepackage{amssymb}
\usepackage{graphicx}
\usepackage{epsf}
\usepackage[T1]{fontenc}
\usepackage[latin9]{inputenc}
\usepackage[letterpaper]{geometry}
\usepackage{setspace}
\usepackage{epstopdf}

\usepackage{rotating}
\makeatletter

\providecommand{\tabularnewline}{\\}

\usepackage{epsf}\@ifundefined{definecolor}
 {\usepackage{color}}{}

\makeatletter

\providecommand{\tabularnewline}{\\}

\usepackage{epsf}

\newcommand{\captionfonts}{\footnotesize}\makeatletter
\long\def\@makecaption#1#2{%
  \vskip\abovecaptionskip
  \sbox\@tempboxa{{\captionfonts #1: #2}}%
  \ifdim \wd\@tempboxa >\hsize
    {\captionfonts #1: #2\par}
  \else
    \hbox to\hsize{\hfil\box\@tempboxa\hfil}%
  \fi
  \vskip\belowcaptionskip}\makeatother

\newtheorem{theorem}{Theorem}

\newtheorem{lemma}{Lemma}

\newtheorem{ass}{Assumption}

\begin{document}

\centerline{\large{Conditional Inference with a Functional Nuisance Parameter}}

\centerline{By Isaiah Andrews\footnote{Harvard Society of Fellows.  Harvard Department of Economics, Littauer Center M39, Cambridge, MA 02138. Email iandrews@fas.harvard.edu.  NSF Graduate Research Fellowship support under grant number 1122374 is gratefully acknowledged.}
 and Anna Mikusheva \footnote{
Department of Economics, M.I.T., 77 Massachusetts Avenue, E18-224, Cambridge, MA, 02139. Email: amikushe@mit.edu.
      Financial support from the Castle-Krob
      Career Development Chair and the Sloan Research Fellowship is gratefully acknowledged.
      We thank Alex Belloni, Victor Chernozhukov, Kirill Evdokimov, and Martin Spindler for helpful discussions.
}}

\centerline{Abstract}

 {\small{This paper shows that the problem of testing hypotheses in moment condition models without any assumptions about identification may be considered
as a problem of testing with an infinite-dimensional nuisance parameter.  We introduce a sufficient statistic for this nuisance parameter and propose conditional tests.  These conditional tests have uniformly correct asymptotic size for a large class of models and test statistics.  We apply our approach to construct tests based on quasi-likelihood ratio statistics, which we show are efficient in strongly identified models and perform well relative to existing alternatives in two examples.
}\\
\texttt{Key words: weak identification, similar test, conditional inferences} }

\begin{center}

\end{center}

\section{Introduction}\label{section-introduction}
Many econometric techniques identify and draw inferences about a structural parameter $\theta$ based on a set of moment equalities.  In particular, many models imply that some function of the data and model parameters has mean zero when evaluated at the true parameter value $\theta_0$.  The current econometric literature devotes a great deal of energy to investigating whether a given set of moment restrictions suffices to uniquely identify the parameter $\theta$, and to studying inference under different identification assumptions.  The goal of this paper is to develop techniques for testing that a specific value $\theta_0$ is consistent with the data using a wide variety of test statistics, without making any assumption about the point identification or strength of identification of the model.

We treat moment equality models as having a functional nuisance parameter.  Much work in econometrics focuses on $\theta$ as the unknown model parameter, typically belonging to a finite-dimensional parameter space.  This is consistent with the tradition from classical statistics, which studied fully-parametric models where the unknown parameter $\theta$ fully described the distribution of the data.  By contrast, in moment condition models the joint distribution of the data is typically only partially specified, and in particular the mean of the moment condition at values $\theta$ other than $\theta_0$ is typically unknown. In light of this fact we suggest re-considering the parameter space in these semi-parametric models, and view the mean function as an unknown (and often infinite-dimensional) parameter. The structural parameter $\theta_0$ corresponds to a zero of this unknown function, and any hypothesis about $\theta_0$ can be viewed as a composite hypothesis with an infinite-dimensional nuisance parameter, specifically the value of the mean function for all other values $\theta$.  The mean function determines the identification status of the structural parameter $\theta$, thus treating the mean function as a parameter allows us to avoid making assumptions about identification.
Corresponding to this infinite-dimensional parameter, we base inference on observation of an infinite-dimensional object, namely the stochastic process given by the sample moment function evaluated at different parameter values $\theta$.

This perspective allows us to study the behavior of a wide variety of test statistics for the hypothesis that the mean function is equal to zero at $\theta_0$.  In a point-identified setting this hypothesis corresponds to testing that $\theta_0$ is the true parameter value, while when point identification fails it corresponds to testing that $\theta_0$ belongs to the identified set.  The existing literature proposes a number of tests for this hypothesis but most of these procedures
depend on the observed process only through its value, and potentially derivative, at the point $\theta_0$.  Examples include the Anderson-Rubin statistic, Kleibergen (2005)'s K statistic, and generalizations and combinations of these.  A major reason for restricting attention to statistics which depend only on behavior local to $\theta_0$ is that the distribution of these statistics is independent of the unknown mean function, or depends on it only through a finite-dimensional parameter.  Unfortunately, however, restricting attention to the behavior of the process local to $\theta_0$ ignores a great deal of information and so may come at a significant cost in terms of power.  Further, this restriction rules out many test statistics known to have desirable power properties in other settings.  In contrast to the previous literature, our approach allows us to consider test statistics which depend on the full path of the observed process.

To construct tests based on these statistics, we introduce a sufficient statistic for the unknown mean function and condition inference on the realization of this sufficient statistic.  The idea of conditioning on a sufficient statistic for a nuisance parameter is a longstanding tradition in statistics and was popularized in econometrics by Moreira (2003), which applied this idea in weakly-identified linear instrumental variables models.  The contribution of this paper is to show how this technique may be applied in contexts with an infinite-dimensional nuisance parameter, allowing its use in a wide range of econometric models.  Since the nuisance parameter in our context is a function, our sufficient statistic is a stochastic process. Our proposed approach to testing is computationally feasible and is of similar difficulty as other simulation-based techniques such as the bootstrap.

One statistic allowed by our approach is the quasi-likelihood ratio (QLR) statistic.  This statistic makes use of the full path of the observed stochastic process and its distribution under the null depends on the unknown mean function, which greatly limited its use in the previous literature on inference with nonstandard identification.  At the same time, one may expect QLR tests to have desirable power properties: in well identified (point identified and strongly identified) models QLR tests are asymptotically efficient, while they avoid the power deficiencies of Kleibergen (2005)'s K and related tests under weak identification.  Moreover, in linear IV with homoskedastic errors Andrews, Moreira, and Stock (2006) showed that Moreira (2003)'s conditional likelihood ratio (CLR) test, which corresponds to the conditional QLR test in that context, is nearly uniformly most powerful in an important class of tests.

Conditioning on a sufficient statistic for a nuisance parameter, while widely applied, may incur loss of power by restricting the class of tests permitted.  We show, however, that no power loss is incurred in well identified models as in this case our conditional QLR test is asymptotically equivalent to the unconditional QLR test and thus is efficient.  We also point out that if one is interested in similar tests (that is, tests with exactly correct size regardless of the mean function) and the set of mean functions is rich enough, all similar tests are conditional tests of the form we consider.

To justify our approach we show that for a large class of test statistics the conditional tests we propose have uniformly correct asymptotic size over a broad class of models which imposes no restriction on the mean function, and so includes a wide range of identification settings.  We further extend these results to allow for concentrating out well-identified structural nuisance parameters.

We apply our approach to inference on the coefficients on the endogenous regressors in the quantile IV model studied by Chernozhukov and Hansen (2005, 2006, 2008) and Jun (2008).  We examine the performance of the conditional QLR test in this context and find that it has desirable power properties relative to alternative approaches.  In particular, unlike Anderson-Rubin-type tests the conditional QLR test is efficient under strong identification, while unlike tests based on the K statistic it does not suffer from non-monotonic power under weak identification.

As an empirical application of our method, we compute confidence sets for nonlinear Euler Equation parameters based on US data. We find that our approach yields much smaller confidence sets than existing alternatives, and in particular allows us to rule out high values of risk aversion allowed by alternative methods.

In Section \ref{section-setting} we introduce our model and discuss the benefits of formulating the problem using an infinite-dimensional nuisance parameter. Section \ref{section- conditional approach} explains and justifies our conditioning approach and relates our results to previous work. Section \ref{section-asymptotics} establishes the uniform asymptotic validity of our method and proves the asymptotic efficiency of the conditional QLR test in strongly identified settings, while Section \ref{section- concentrating out} discusses the possibility of concentrating out well-identified nuisance parameters. Section \ref{section- numerical} reports simulations on the power properties of the conditional QLR test in a quantile IV model and gives confidence sets for nonlinear Euler equation parameters based on US data, and Section \ref{section: conclude} concludes. Some proofs and additional results may be found in a Supplementary Appendix available on the authors' web-sites.

In the remainder of the paper we denote by $\lambda_{\min}\left(A\right)$ and $\lambda_{\max}\left(A\right)$
the minimal and maximal eigenvalues of a square matrix $A$, respectively, while $\|A\|$ is the operator norm for a matrix and the Euclidean norm for a vector.

\section{Models with functional nuisance parameters}\label{section-setting}

Many testing problems in econometrics can be recast as tests that a vector-valued random function of model parameters has mean zero at a particular point.
Following Hansen (1982) suppose we have an economic model which implies
that some $k\times 1$-dimensional function $\varphi\left(X_{t};\theta\right)$ of the data and the $q\times 1$-dimensional parameter $\theta$
has mean zero when evaluated at the true parameter value $\theta_0$,
$E\left[\varphi\left(X_{t},\theta_0\right)\right]=0.$ Define $g_{T}(\cdot)=\frac{1}{\sqrt{T}}\sum_{t=1}^{T}\varphi\left(X_{t},\cdot\right)$ and let
$m_{T}(\cdot)=E\left[g_{T}\left(X_{t},\cdot\right)\right].$  Under mild conditions (see e.g. Van der Vaart and Wellner (1996)), empirical process theory implies that
\begin{equation}
g_{T}(\theta)=m_T(\theta)+G(\theta)+r_{T}(\theta),\label{eq:Gaussian Approximation}
\end{equation}
where  $G(\cdot)$ is a  mean-zero Gaussian process with consistently estimable
covariance function $\Sigma(\theta,\tilde\theta)=EG(\theta)G(\tilde\theta)'$,    and $r_{T}$
is a residual term which is uniformly negligible for large $T$.  We are interested in  testing that $\theta_0$ belongs to the identified set, which is equivalent to testing $H_{0}:m_T\left(\theta_{0}\right)=0$, without any assumption on identification of the parameter $\theta$.

This paper considers (\ref{eq:Gaussian Approximation}) as a model with an infinite-dimensional nuisance parameter, namely $m_T(\theta)$ for $\theta\neq\theta_0$.  Thus our perspective differs from the more classical approach which focuses on $\theta$ as the model parameter.  This classical approach may be partially derived from the use of parametric models in which $\theta$ fully specifies the distribution of the data.  By contrast many of the methods used in modern econometrics, including GMM, only partially specify the distribution of the data, and the behavior of $m_T(\theta)$ for $\theta$ outside of the identified set is typically neither known nor consistently estimable. To formally describe the parameter space for $m_T$, let $\mathcal{M}$ be the set of functions $m_T(\cdot)$ that may arise in a given model, and let $\mathcal{M}_0$ be the subset of $\mathcal{M}$ containing those functions satisfying $m_T(\theta_0)=0$.  The hypothesis of interest may be formulated as $H_0:m_T\in\mathcal{M}_0$, which is in general a composite hypothesis with a non-parametric nuisance parameter.

The distribution of most test statistics under the null depends crucially on the nuisance function $m_T(\cdot)$. For example the distribution of quasi-likelihood ratio (QLR) statistics, which for $\widehat\Sigma$ an estimator of $\Sigma$ takes the form
\begin{align}\label{eq: defn of qlr}
QLR=g_T(\theta_0)'\widehat\Sigma(\theta_0,\theta_0)^{-1}g_T(\theta_0)- \inf_\theta g_T(\theta)'\widehat\Sigma(\theta,\theta)^{-1}g_T(\theta),
\end{align}
depends in complex ways on the true unknown function $m_T(\cdot)$, except in special cases like the strong identification assumptions introduced in Section \ref{section- strong id}. The same is true of Wald- or t-statistics, or of statistics analogous to QLR constructed using a weighting other then $\widehat\Sigma(\theta,\theta)^{-1}$, which we call QLR-type statistics.
In the literature to date the dependence on $m_T$ has greatly constrained the use of these statistics in non-standard settings, since outside of special cases (for example linear IV, or the models studied by Andrews and Cheng (2012)) there has been no way to calculate valid critical values.

Despite these challenges there are a number of tests in the literature that control size for all values of the infinite-dimensional nuisance parameter $m_T(\cdot)$. One well-known example is the S-test of Stock and Wright (2000), which is based on the statistic $S=g_T(\theta_0)'\widehat\Sigma(\theta_0,\theta_0)^{-1}g_T(\theta_0)$.  This statistic is a generalization of the Anderson-Rubin statistic and is asymptotically $\chi^2_k$ distributed for all $m_T\in\mathcal{M}_0$.  Other examples include Kleibergen (2005)'s K test and its generalizations.  Unfortunately, these tests often have deficient power in over-identified settings or when identification is weak, respectively.  Several authors have also suggested statistics intended to mimic the behavior of QLR in particular settings, for example the GMM-M statistic of Kleibergen (2005), but the behavior of these statistics differs greatly from true QLR statistics in many contexts of interest.

\paragraph{Example 1.}
Consider the nonlinear Euler equations studied by Hansen and Singleton (1982).
The moment function identifying the discount factor $\delta$ and the coefficient of relative risk-aversion $\gamma$ is
$$
g_T(\theta)=\frac{1}{\sqrt{T}}\sum_{t=1}^T\left( \delta\left(\frac{C_{t}}{C_{t-1}}\right)^{-\gamma}R_{t}-1\right)Z_{t}, ~~~~
\theta=\left(\delta,\gamma\right),
$$
where $C_{t}$ is consumption in period $t$, $R_{t}$ is an asset return
from period $t-1$ to $t,$ and $Z_{t}$ is a vector of instruments
measurable with respect to information at $t-1$. Under moment and mixing conditions (see for example Theorem 5.2 in Dedecker and Louhici (2002)), the demeaned process $g_T(\cdot)-Eg_T(\cdot)$ will converge uniformly to a Gaussian process.

For true parameter value $\theta_0=\left(\delta_0,\gamma_0\right)$ we have $m_T(\theta_0)=Eg_T(\theta_0)=0$.
The value of $m_T(\theta)=Eg_T(\theta)$ for $\theta\neq\theta_0$ is in general unknown and depends in a complicated way on the joint distribution of the data, which is typically neither known nor explicitly modeled.
Further, $m_T(\theta)$ cannot be consistently estimated. Consequently the distribution of QLR and many other statistics which depend on $m_T(\cdot)$ are unavailable unless one is willing to assume the model is well-identified, which is contrary to extensive evidence suggesting identification problems in this context.
 $\Box$

\subsection{The mean function $m_T$ in examples}\label{section-examples}
Different econometric settings give rise to different mean functions $m_T(\cdot)$, which in turn determine the identification status of $\theta$.  In set-identified models the identified set $\{\theta:m_T(\theta)=0\}$ might be a collection of isolated points or sets, or even the whole parameter space.  In well-identified settings, by contrast, $m_T(\cdot)$ has a unique zero and increases rapidly as we move away from this point, especially as $T$ becomes large.  Common models of weak identification imply that even for $T$ large $m_T(\cdot)$ remains bounded over some non-trivial region of the parameter space.

Consider for example the classical situation (as in Hansen (1982)) where the function $E\varphi(X_t,\cdot)$ is fixed and continuously differentiable with a unique zero at $\theta_0$, and the Jacobian $\frac{\partial E\varphi(X_t,\theta_0)}{\partial\theta}$ has full rank. This is often called a \emph{strongly identified case}, and (under regularity conditions) will imply the strong identification assumptions we introduce in Section \ref{section- strong id}.
In this setting the function $m_T(\theta)=\sqrt{T}E\varphi(X_t,\theta)$ diverges to infinity outside of $1/\sqrt{T}$ neighborhoods of $\theta_0$ as the sample size grows. Many  statistics, like Wald or QLR-type statistics, use $g_T(\cdot)$ evaluated  only at some estimated value $\widehat\theta$ and $\theta_0$, and thus in the classical case they depend on $g_T$ only through its behavior on a $1/\sqrt{T}$ neighborhood of the true $\theta_0$. Over such neighborhoods
$m_T(\cdot)$ is well approximated by  $\sqrt{T}\frac{\partial E\varphi(X_t,\theta_0)}{\partial\theta}(\theta-\theta_0)$, the only unknown component of which, $\frac{\partial E\varphi(X_t,\theta_0)}{\partial\theta},$ is usually consistently estimable. Reasoning along these lines, which we explore in greater detail in Section \ref{section- strong id}, establishes the asymptotic validity of classical tests under strong identification.
Thus in strongly identified models the nuisance parameter problem we study here does not arise.

In contrast to the strongly-identified case, \emph{weakly identified} models are often understood as those in which even for $T$ large the mean function $m_T$ fails to dominate the Gaussian process $G$ over a substantial part of the parameter space.
Stock and Wright (2000) modeled this phenomenon using a drifting sequence
of functions. In particular, a simple case of the Stock and Wright (2000) embedding indexes the data-generating process by the sample size
 and assumes that while the variance of the moment condition is asymptotically constant, the expectation of the moment condition shrinks at rate $1/\sqrt{T}$, so $E\varphi(X_t,\theta)=E_T\varphi(X_t,\theta)=\frac{1}{\sqrt{T}}f(\theta)$ for a fixed function $f(\theta)$. In this case $m_T(\theta)=f(\theta)$ is  unknown and cannot be consistently estimated, consistent estimation of $\theta_0$ is likewise impossible, and the whole function $m_T(\cdot)$ is important for the distribution of QLR-type statistics.

By treating $m_T$ as a nuisance parameter, our approach avoids making any assumption on its behavior.  Thus, we can treat
both the strongly-identified case described above and the weakly-identified sequences studied by Stock and Wright (2000), as well as set identified models and a wide array of other cases.  As we illustrate below this is potentially quite important, as the set $\mathcal{M}$ of mean functions can be extremely rich in examples.

We next  discuss the sets $\mathcal{M}$ in several examples. As a starting point we consider the linear IV model, where the nuisance function can be reduced to a finite-dimensional vector of nuisance parameters, and then consider examples with genuine functional nuisance parameters.

\paragraph{Example 2. (Linear IV)}
Consider a linear IV model where the data consists of i.i.d. observations on an outcome variable
$Y_{t}$, an endogenous regressor $D_{t}$,
and a vector of instruments $Z_{t}$. Assume that the identifying moment condition is $E\left[\left(Y_{t}-D_{t}'\theta_0\right)Z_{t}\right]=0$.
This implies that $m_{T}\left(\theta\right)=\sqrt{T}E\left[Z_{t}D_{t}'\right]\left(\theta_0-\theta\right)$ is a linear function.
If $E\left[Z_{t}D_{t}'\right]$ is a fixed matrix of full column rank, then $\theta_0$
is point identified and can be consistently estimated using two-stage-least-squares, while if $E\left[Z_{t}D_{t}'\right]$ is
of reduced rank the identified set is a hyperplane of dimension
equal to the rank deficiency of $E\left[Z_{t}D_{t}'\right]$.
Staiger and Stock (1997) modeled weak instruments by considering
a sequence of data-generating processes such that $E\left[Z_{t}D_{t}'\right]=\frac{C}{\sqrt{T}}$
for a constant unknown matrix $C$. Under these sequences the function $m_T(\theta)=C(\theta_0-\theta)$ is linear and governed by the unknown (and not consistently estimable) parameter $C$. $\Box$

In contrast to the finite-dimensional nuisance parameter obtained in linear IV, in nonlinear models the space of nuisance parameters $m_T(\cdot)$ is typically of infinite dimension.

\paragraph{Example 1 (continued).}
In the Euler equation example discussed above,
$$
m_{T}\left(\theta\right)=\sqrt{T}E\left[\left(\delta\left(1+R_{t}\right) \left(\frac{C_{t}}{C_{t-1}}\right)^{-\gamma}- 1\right)Z_{t}\right].
$$
Assume for a moment that $\delta$ is fixed and known and  that $R_{t}$ and
$Z_{t}$ are constant. In this simplified case the function $m_{T}(\gamma)$ is a linear transformation of
the moment generating function of $\log\left(C_{t}/C_{t-1}\right)$,
implying that the set $\mathcal{M}_0$ of mean functions is at least
as rich as the set of possible distributions for consumption growth consistent with the null. $\Box$

\paragraph{Example 3.}
In a  nonlinear IV models with the moment condition $$E\left[\left(Y_{t}-f\left(D_{t},\theta\right)\right)Z_{t}\right]=0$$
the mean function has the form
$$
m_{T}\left(\theta\right)=\sqrt{T}E\left[\left(f\left(D_{t},\theta_0\right)-f\left(D_{t},\theta\right)\right)Z_{t}\right].
$$
The set of nuisance parameters $\mathcal{M}_0$ will in general depend on the structure of the function
$f$. For example, if $f$ is multiplicatively separable in data
and parameters, so $f\left(D_{t},\theta\right)=f_{1}\left(D_{t}\right)'f_{2}\left(\theta\right),$
then we can write $m_{T}\left(\theta\right)=\sqrt{T}E\left[Z_{t}f_{1}\left(D_{t}\right)'\right](f_{2}\left(\theta_0\right)-f_{2}\left(\theta\right)),$
and similar to the linear IV model the moment function will be governed
by the finite-dimensional nuisance parameter $\sqrt{T}E\left[Z_{t}f_{1}\left(D_{t},C_{t}\right)'\right]$.
In more general models, however, the function $m_{T}(\cdot)$ may depend
on the distribution of the data in much richer ways, leaving us with
an infinite-dimensional nuisance parameter.  $\Box$

 Our results
also apply outside the GMM context so long as one has a model described by (\ref{eq:Gaussian Approximation}).
In Section \ref{section- quantile IV}, for example, we apply our results to a quantile IV
where we plug in estimates for nuisance parameters. Our results can likewise be applied to
 the simulation-based moment conditions considered
in McFadden (1989), Pakes and Pollard (1989), and the subsequent literature.
More recently Schennach (2014) has shown that models with latent variables
can be expressed using simulation-based moment conditions, allowing
the treatment of an enormous array of additional examples including
game-theoretic, moment-inequality, and measurement-error models within
the framework studied in this paper.

\section {Conditional approach}\label{section- conditional approach}

To construct tests we introduce a sufficient statistic for $m_T(\cdot)\in\mathcal{M}_0$
and suggest conditioning inference on this statistic, thereby eliminating dependence on the nuisance parameter. Moreira (2003) showed that the conditioning approach could
be fruitfully applied to inference in linear instrumental variables
models, while Kleibergen (2005)  extended this approach to GMM
statistics which depend only on $g_T(\cdot)$ and its derivative
both evaluated at $\theta_0$. In this section we show that conditional tests can be applied far more broadly.  We first introduce our approach and describe how to calculate critical values,
 then justify our procedure in a limit problem.
 In Section \ref{section-asymptotics} we show that our tests are uniformly asymptotically  correct under more general assumptions.

\subsection{Conditional inference}

Consider model (\ref{eq:Gaussian Approximation}), and let $\widehat\Sigma(\cdot,\cdot)$ be a consistent estimator of covariance function $\Sigma(\cdot,\cdot)$.
Let us introduce the process
\begin{align}\label{eq: defn of h}
h_T\left(\theta\right)=H(g_T,\widehat\Sigma)(\theta)=g_T\left(\theta\right)- \widehat\Sigma\left(\theta,\theta_{0}\right) \widehat\Sigma\left(\theta_{0},\theta_{0}\right)^{-1}g_T\left(\theta_{0}\right).
\end{align}
We show in Section \ref{section- limit problem} that this process is a sufficient statistic for $m_T(\cdot)\in\mathcal{M}_0$ in the limit problem
where the residual term in (\ref{eq:Gaussian Approximation}) is exactly zero and the covariance of $G(\cdot)$ is  known ($\widehat\Sigma(\cdot,\cdot)=\Sigma(\cdot,\cdot)$).
Thus the conditional distribution of any test statistic $R=R(g_T,\widehat\Sigma)$ given $h_T(\cdot)$ does not depend on the nuisance parameter $m_T(\cdot)$. Following the classical conditioning approach (see e.g. Lehmann and Romano (2005)) we create a test based on statistic $R$ by pairing it with conditional critical values that depend on the process $h_T(\cdot)$.

To simulate the conditional distribution of statistic $R$ given $h_T(\cdot)$ we take independent draws $\xi^*\sim N(0,\widehat\Sigma(\theta_0,\theta_0))$ and produce simulated processes
\begin{align}\label{eq: simulating quantiles}
g^*_T\left(\theta\right)=h_T\left(\theta\right)+\widehat\Sigma\left(\theta,\theta_{0}\right) \widehat\Sigma\left(\theta_{0},\theta_{0}\right)^{-1}\xi^*.
\end{align}
We then calculate  $R^*=R(g_T^*,\widehat\Sigma)$, which represents a random draw from the conditional distribution of $R$ given $h_T$ under the null (in the limit problem). To calculate the conditional $(1-\alpha)$-quantile of $R$ to use as a critical value, we can thus simply take the $(1-\alpha)$-quantile of $R^*$, which is straightforward to approximate by simulation.

\subsection{Limit problem}\label{section- limit problem}

In this section we consider a limit problem that abstracts from some finite-sample features but leaves the central challenge of inference with an infinite-dimensional nuisance parameter intact.
Consider a statistical experiment in which we observe the process $g_T(\theta)=m_T(\theta)+G(\theta),$
where $m_T(\cdot)\in\mathcal{M}$ is an unknown deterministic mean function, and $G(\cdot)$ is a mean-zero Gaussian process with known covariance $\Sigma(\theta,\tilde\theta)=EG(\theta)G(\tilde\theta)'$. We again assume that $\mathcal{M}$ is the set of potential mean functions, which is in general infinite-dimensional, and wish to test the hypothesis $H_0:m_T(\theta_0)=0$.

Lemma \ref{lemma- sufficient statistic} below shows that the process $h_T(\cdot)$ is a sufficient statistic for the unknown function $m_T(\cdot)$ under the null $m_T(\cdot)\in\mathcal{M}_0$. The validity of this statement hinges on the observation that under the null the process $g_T(\cdot)$ can be decomposed into two independent, random components- the process $h_T(\cdot)$  and the random vector $g_T(\theta_0)$:
\begin{align}\label{eq: g as a function of h}
g_T(\theta)=h_T\left(\theta\right)+\Sigma\left(\theta,\theta_{0}\right) \Sigma\left(\theta_{0},\theta_{0}\right)^{-1}g_T(\theta_0),
\end{align}
with the important property that the distribution of $g_T(\theta_0)\sim N(0,\Sigma(\theta_0,\theta_0))$ does not depend on the nuisance parameter $m_T(\cdot)$.
In particular, this implies that the conditional distribution of any functional of $g_T(\cdot)$ given $h_T(\cdot)$ does not depend on $m_T(\cdot)$.

Assume we wish to construct a test that rejects the null hypothesis when the statistic $R=R(g_T,\Sigma)$, calculated using the observed $g_T(\cdot)$ and the  known covariance $\Sigma(\cdot,\cdot)$, is large.  Define the conditional critical value function $c_{\alpha}\left(h_T\right)$
by
$$
c_{\alpha}\left(\tilde{h}\right)=\min\left\{ c:P\left\{R\left(g_T,\Sigma\right)>c\mid h_T=\tilde{h}\right\}\leq\alpha \right\}.
$$
Note that the conditional quantile $c_\alpha(\cdot)$ does not depend on the unknown $m_T(\cdot)$, and that for any realization of $h_T(\cdot)$ it can be easily simulated as described above.

\begin{lemma}\label{lemma- sufficient statistic}
In the limit problem the test that rejects the null hypothesis $H_{0}:m_T\in\mathcal{M}_{0}$ when $R(g_T,\Sigma)$ exceeds the random critical value $c_\alpha(h_T)$  has correct size. If the conditional distribution of $R$  given $h_T$ is continuous almost surely  then the test is conditionally similar given $h_T(\cdot)$. In particular, in this case
for any $m_T\in\mathcal{M}_0$ we have that almost surely
$$P\left\{R(g_T,\Sigma)>c_{\alpha}(h_T)\mid h_T(\cdot)\right\}=P\left\{R(g_T,\Sigma)>c_{\alpha}(h_T)\right\}=\alpha.$$

\end{lemma}

The critical value $c_\alpha(h_T)$ is a random variable, as it depends on random process $h_T$. Under an almost sure continuity assumption the proposed test is conditionally similar, in that is has conditional size $\alpha$ for almost every realization of $h_T$.

Conditional similarity is a very strong restriction and may be hard to justify in some cases as it greatly reduces the class of possible tests.  If, however, one is interested in similar tests (tests with exact size $\alpha$ regardless of the value of the nuisance parameter), all such tests will automatically be conditionally similar given a sufficient statistic if the family of distributions for the sufficient statistic under the null is boundedly complete- we refer the interested reader to Lehmann and Romano (2005) and Moreira (2003) for further discussion of this point.

If the parameter space for $\theta$ is finite ($\Theta=\{\theta_0,\theta_1,...,\theta_n\}$) the conditions for bounded completeness are well-known and easy to check.  In particular, in this case our problem reduces to that of observing a $k(n+1)$-dimensional Gaussian vector $g_T=(g_T(\theta_0)',...,g_T(\theta_n)')'$ with unknown mean $(0,\mu_1'=m_T(\theta_1)',...,\mu_n'=m_T(\theta_n)')'$  and known covariance.  If the set $\mathcal{M}$ of possible values for the nuisance parameter $(\mu_1',...,\mu_n')'$ contains a rectangle with a non-empty interior then the family of distributions for $h_T$ under the null is boundedly complete, and all similar tests are conditionally similar given $h_T$.  A generalization of this statement to cases with infinite-dimensional nuisance parameters is provided in the Supplementary Appendix.

While similarity is still a strong restriction, similar tests have been shown to perform well in other weakly identified contexts, particularly in linear IV: see Andrews Moreira and Stock (2008).  On a practical level, as we detail below the presence of the infinite-dimensional nuisance parameter $m_T\in\mathcal{M}_0$ renders many other approaches to constructing valid tests unappealing in the present context, as alternative approaches greatly restrict the set of models considered, the set of test statistics permitted, or both.

\subsection{Relation to the literature}

Moreira (2003) pioneered the conditional testing approach in  linear IV models with homoskedastic errors, which are a special case of our Example 2.
If we augment Example 2 by assuming that the instruments $Z_t$ are non-random and the reduced form errors are Gaussian with mean zero and known covariance matrix $\Omega$, we obtain a model satisfying the assumptions of the limit problem in each sample size.  In particular, for each $T$ we observe the process
$g_T(\theta)=\frac{1}{\sqrt{T}}\sum_{t=1}^T(Y_t-D_t'\theta)Z_t$, which is Gaussian with mean function $m_T(\theta)=\frac{1}{\sqrt{T}}\sum_{t=1}^TE[Z_tD_t'](\theta_0-\theta)$ and covariance function
$$
\Sigma(\theta,\tilde\theta)=\left(\frac{1}{T}\sum_{t=1}^{T}Z_tZ_t'\right)(1,-\theta)\Omega(1,-\tilde\theta)'.
$$
In this case both the mean function $m_T(\cdot)$ and the process $g_T(\cdot)$ are linear, and so belong to a finite-dimensional space.  The process $h_T(\cdot)$ is likewise linear
in this model, and its coefficient of linearity is proportional to the statistic that Moreira (2003) called $T$ and used as the basis of his conditioning technique.
Thus, the conditioning we propose is equivalent to that suggested by Moreira (2003) in linear IV, and our approach is a direct generalization of Moreira (2003) to nonlinear models.  Consequently, when applied to the QLR statistic in homoskedastic linear IV, our approach yields the CLR test of Moreira (2003), which Andrews, Moreira, and Stock (2006) shows is nearly a uniformly most powerful test in a class of invariant similar two-sided tests in the homoskedatic Gaussian linear IV model.

Kleibergen (2005) generalized the conditioning approach of Moreira (2003) to some statistics for potentially nonlinear GMM models. Kleibergen (2005) restricts attention to statistics which depend on the data only through  $g_T(\theta_0)$ and $\frac{d}{d\theta}g_T(\theta_0)$, which he assumes to be jointly Gaussian in the limit experiment.  To produce valid tests he pairs these statistics with critical values calculated by conditioning on a statistic he called $D_T$, which can be interpreted as the part of  $\frac{d}{d\theta}g_T(\theta_0)$ which is independent of $g_T(\theta_0)$.  One can easily show, however, that in the limit problem Kleibergen's $D_T$ is the negative of $\frac{d}{d\theta}h_T(\theta_0)$.  Moreover, one can decompose $h_T(\cdot)$ into the random matrix $\frac{d}{d\theta}h_T(\theta_0)$ and a process which is independent of both $\frac{d}{d\theta}h_T(\theta_0)$ and $g_T(\theta_0)$, so the conditional distribution of any function of $g_T(\theta_0)$ and $\frac{d}{d\theta}g_T(\theta_0)$ given $h_T(\cdot)$ is simply its conditional distribution given $\frac{d}{d\theta}h_T(\theta_0)$.  Thus, for the class of statistics considered in Kleibergen (2005) our conditioning approach coincides with his.  Unlike Kleibergen (2005), however, our approach can treat statistics which depend on the full process $g_T(\cdot)$, not just on its behavior local to the null.  In particular our approach allows us to consider QLR statistics, which are outside the scope of Kleibergen's approach in nonlinear models.  Kleibergen (2005) introduces what he terms a GMM-M statistic, which coincides with the CLR statistic in homoskedastic linear IV and is intended to extend the properties of the CLR statistic to more general settings, but this statistic unfortunately has behavior quite different from a true QLR statistic in some empirically relevant settings, as we demonstrate in an empirical application to the Euler equation example 1 in Section \ref{subsection- Euler Equation}.

\paragraph{Unconditional tests with nuisance parameters.} In models with finite-dimensional nuisance parameters, working alternatives to the conditioning approach include least favorable and Bonferroni critical values. Least favorable critical values search over the space of nuisance parameters to maximize the ($1-\alpha$)-quantile of the test statistic, and this approach was successfully implemented by Andrews and Guggenberger (2009) in models with a finite-dimensional nuisance parameter.  Unfortunately, however, in cases with a functional nuisance parameter the least-favorable value is typically unknown and a simulation search is computationally infeasible, rendering this approach unattractive.  Bonferroni critical values are similar to least favorable ones, save that instead of searching over whole space of nuisance parameters we instead search only over some preliminary confidence set.  Again, absent additional structure this approach is typically only feasible when the nuisance parameter is of finite dimension.  Relatedly, Andrews and Cheng (2012) show that in the settings they consider the behavior of estimators and test statistics local to a point of identification failure are controlled by a finite-dimensional nuisance parameter and use this fact to construct critical values for QLR and Wald statistics which control size regardless of the value of this parameter.

Common ways to calculate critical values in other contexts include subsampling and the bootstrap.  Both of these approaches are known to fail to control size for many test statistics even in cases with finite-dimensional nuisance parameters, however (see Andrews and Guggenberger (2010)), and thus cannot be relied on in the present setting.  Indeed, it is straightforward to construct examples demonstrating that neither subsampling nor the bootstrap yields valid critical values for the QLR statistic in general.

\section{ Asymptotic behavior of conditional tests}\label{section-asymptotics}

\subsection{ Uniform validity}
The limit problem studied in the previous section assumes away many finite-sample features relevant in empirical work, including non-Gaussianity of $g_T$ and error in estimating the covariance function $\Sigma$.  In this section we extend our results to allow for these issues, and show that our conditioning approach yields uniformly asymptotically valid tests over large classes of models in which the observed process $g_T(\cdot)$ is uniformly asymptotically Gaussian.

Let $P$ be a probability measure describing the distribution of $g_T(\cdot)$, where $T$ denotes the sample size. For each probability law $P$ there is a deterministic mean function $m_{T,P}(\cdot)$, which will in many cases be the expectation $E_Pg_T(\cdot)$ of the process $g_T(\cdot)$ under $P$. We assume that the difference $g_T(\cdot)-m_{T,P}(\cdot)$ converges to a mean zero Gaussian process $G_P(\cdot)$ with covariance function $\Sigma_P(\cdot,\cdot)$ uniformly over the family $\mathcal{P}_0$ of distributions consistent with the null.  We formulate this assumption using bounded Lipshitz convergence- see Van der Vaart and Wellner (1996) for the equivalence between bounded Lipshitz convergence and weak convergence of stochastic processes.  For simplicity of notation we suppress the subscript $P$ in all expressions:

\begin{ass} \label{Ass: Uniform CLT} The difference $g_{T}\left(\cdot\right)-m_T\left(\cdot\right)$
converges  to a Gaussian process $G(\cdot)$ with mean zero and covariance function $\Sigma(\cdot,\cdot)$ uniformly over $P\in\mathcal{P}_0$, that is:
\[
\lim_{T\to\infty}\sup_{P\in\mathcal{P}_{0}} \sup_{f\in BL_{1}} \left\|E\left[f\left(g_{T}-m_T\right)\right] -E\left[f\left(G\right)\right]\right\|=0,
\]
where $BL_1$ is the set of functionals with Lipshitz constant and supremum norm bounded above by one.
\end{ass}

\begin{ass} \label{Ass: Bounded Covariance} The covariance function $\Sigma(\cdot,\cdot)$ is uniformly bounded and positive definite:
\[
1/\bar{\lambda}\leq\inf_{P\in\mathcal{P}_{0}}\inf_{\theta\in\Theta}\lambda_{\min}\left(\Sigma \left(\theta,\theta\right)\right)\leq
\sup_{P\in\mathcal{P}_{0}}\sup_{\theta\in\Theta}\lambda_{\max} \left(\Sigma\left(\theta,\theta\right)\right)\le\bar{\lambda},
\]
for some finite $\bar{\lambda}>0$.
\end{ass}

\begin{ass}\label{Ass: consistency of covariance} There is a uniformly consistent estimator $\widehat\Sigma(\cdot,\cdot)$ of the covariance function, in that for any $\varepsilon>0$
$$
\lim_{T\to\infty}\sup_{P\in\mathcal{P}_{0}}P\left\{\sup_{\theta,\tilde\theta}\left\| \widehat\Sigma(\theta,\tilde\theta)-\Sigma(\theta,\tilde\theta)\right\| >\varepsilon\right\}=0.
$$
\end{ass}

Suppose we are interested in tests that reject for large values a statistic $R$ which depends on the moment function $g_T(\cdot)$ and the estimated covariance $\widehat\Sigma(\cdot,\cdot)$. Consider process $h_T(\cdot)=H(g_T,\widehat\Sigma)$ defined as in (\ref{eq: defn of h}).
Since the transformation from $(g_T(\cdot),\widehat\Sigma(\cdot,\cdot))$ to $(g_T(\theta_0),h_T(\cdot),\widehat\Sigma(\cdot,\cdot))$ is one-to-one, $R$ can be viewed as a functional of $(g_T(\theta_0),h_T(\cdot),\widehat\Sigma(\cdot,\cdot))$. We require that $R$ be sufficiently continuous with respect to
$(g_T(\theta_0),h_T(\cdot),\widehat\Sigma(\cdot,\cdot))$, which allows  QLR and a number of other statistics but rules out Wald statistics in many models:

\begin{ass}\label{Ass: statistic}  The functional $R(\xi,h(\cdot),\Sigma(\cdot,\cdot))$ is defined for all values $\xi\in\mathbb{R}^k$, all $k$-dimensional functions $h$  with the property that $h(\theta_0)=0$, and all covariance functions $\Sigma(\cdot,\cdot)$ satisfying Assumption \ref{Ass: Bounded Covariance}. For any   fixed $C>0$,  $R(\xi,h,\Sigma)$ is bounded and Lipshitz in $\xi$, $h$, and $\Sigma$  over  the set of $(\xi,h(\cdot),\Sigma(\cdot,\cdot))$  with $\xi'\Sigma(\theta_0,\theta_0)^{-1}\xi\leq C$.
\end{ass}

\begin{lemma}\label{lemma- QLR stat lipshitz}
The QLR statistic defined in (\ref{eq: defn of qlr}) satisfies Assumption. \ref{Ass: statistic}
\end{lemma}

To calculate our conditional critical values, given a realization of $h_T$ we simulate independent draws $\xi\sim N(0,\widehat\Sigma(\theta_0,\theta_0))$ and (letting $P^*$ denote the simulation probability) define
$$
c_\alpha(h_T,\widehat\Sigma)=\inf\left\{c:P^*\left\{\xi: R(\xi,h_T(\cdot),\widehat\Sigma(\cdot,\cdot))\leq c\right\}\geq 1-\alpha\right\}.
$$
The test then rejects if $R(g_T(\theta_0),h_T,\widehat\Sigma)>c_\alpha(h_T,\widehat\Sigma)$.

\begin{theorem}\label{Theorem: uniform asymptotics}
Let Assumptions \ref{Ass: Uniform CLT} - \ref{Ass: statistic} hold, then for any $\varepsilon>0$ we have
$$
\lim_{T\to\infty} \sup_{P\in\mathcal{P}_0} P \left\{ R(g_T(\theta_0),h_T,\widehat\Sigma)>c_\alpha(h_T,\widehat\Sigma)+\varepsilon\right\} \leq \alpha.
$$
\end{theorem}

Theorem \ref{Theorem: uniform asymptotics} shows that our conditional critical value (increased by an arbitrarily small amount) results in a test which is uniformly asymptotically valid over the large class of distributions $\mathcal{P}_0$.  The need for the term $\varepsilon$ reflects the possibility that there may be some sequences of distributions in $\mathcal{P}_0$ under which $R$ converges in distribution to a limit which is not continuously distributed.  If we rule out this possibility, for example assuming that the distribution of $R$ is continuous with uniformly bounded density for all $T$ and all $P\in\mathcal{P}_0$, then the conditional test with $\varepsilon=0$ is uniformly asymptotically similar in the sense of Andrews, Cheng and Guggenberger (2011).

\subsection{Strong identification case}\label{section- strong id}

Restricting attention to conditionally similar tests rules out many procedures and so could come at a substantial cost in terms of power.
In this section, we show that restricting attention to conditionally similar tests does not result in loss of power if the data are in fact generated from a
strongly identified model, by which we mean one satisfying conditions given below.  In particular, we establish that under these conditions our conditional
QLR test is equivalent to the classical QLR test using $\chi^2$ critical values and so retains the efficiency properties of the usual QLR test.

\begin{ass}\label{ass: strong-id 1}
For some sequence of numbers $\delta_T$ converging to zero and each $P\in\mathcal{P}_0$, there exists a sequence of matrices $M_T$ such that for any $\varepsilon>0$:
\begin{itemize}
\item[(i)]$
\lim_{T\to\infty}\inf_{P\in\mathcal{P}_0}\inf_{ \|\theta-\theta_0\|>\delta_T}m_T(\theta)'\Sigma(\theta,\theta)^{-1}m_T(\theta)=\infty,
$
\item[(ii)] $
\lim_{T\to\infty}\sup_{P\in\mathcal{P}_0} \sup_{|\theta-\theta_0|\le\delta_T}|m_T(\theta)-M_T(\theta-\theta_0)|=0,
$
\item[(iii)] $\lim_{T\to\infty}\inf_{P\in\mathcal{P}_0} \delta_T^2\lambda_{min}(M_T'\Sigma(\theta_0,\theta_0)^{-1}M_T)=\infty$,
\item[(iv)] $
\lim_{T\to\infty}\sup_{P\in\mathcal{P}_0}\sup_{\|\theta-\theta_0\|\leq \delta_T}\left\|\Sigma(\theta,\theta)-\Sigma(\theta_0,\theta_0)\right\|=0$ and \\ $\lim_{T\to\infty}\sup_{P\in\mathcal{P}_0}\sup_{\|\theta-\theta_0\|\leq \delta_T}\left\|\Sigma(\theta,\theta_0)-\Sigma(\theta_0,\theta_0)\right\|=0,
$
\item[(v)]$
\lim_{T\to\infty}\sup_{P\in\mathcal{P}_0}P\left\{\sup_{\|\theta-\theta_0\|\leq \delta_T}|G(\theta)-G(\theta_0)|>\varepsilon\right\}=0,
$

\item[(vi)]  There exists a constant $C$ such that
$
\sup_{P\in\mathcal{P}_0}P\left\{\sup_{\theta\in\Theta}|G(\theta)|>C\right\}<\varepsilon.
$
\end{itemize}
\end{ass}

\paragraph{Discussion of Assumption \ref{ass: strong-id 1}.}  Assumption \ref{ass: strong-id 1} defines what we mean by strong identification.  Part (i) guarantees that the moment function diverges outside of a shrinking neighborhood of the true parameter value and, together with assumption (vi), implies the existence of consistent estimators.  Part (ii) requires that the unknown mean function $m_T(\theta)$ be linearizable on a neighborhood of $\theta_0$, which plays a key role in establishing the asymptotic normality of estimators.  Part (iii) follows from parts (i) and (ii) if we require $m_T$ to be uniformly continuously differentiable at $\theta_0$, while parts (iv)-(vi) are regularity conditions closely connected to stochastic equicontinuity.
In particular, (iv) requires that the covariance function be continuous at $\theta_0$, while (v) requires that $G$ be equicontinious at $\theta_0$, and (vi) requires that $G$ be bounded almost surely.

Parts (i)-(iii) of Assumption  \ref{ass: strong-id 1} are straightforward to verify in a classical GMM setting.  Consider a GMM model as in
Section \ref{section-examples} which satisfies Assumptions \ref{Ass: Uniform CLT}-\ref{Ass: Bounded Covariance} with mean function $Eg_T(\theta)=m_T(\theta)=T^{1/2-\alpha}m(\theta)$, where $0\leq\alpha<1/2$ and $m(\theta)$ is a fixed, twice-continuously-differentiable function with $m(\theta)=0$ iff $\theta=\theta_0$. Assume further that  $m(\theta)$ is continuously differentiable at $\theta_0$ with full-rank Jacobian $\frac{\partial}{\partial\theta}m(\theta_0)=M$, and that the parameter space $\Theta$ is compact.  For $\delta_T=T^{-\gamma}$, $\inf_{ \|\theta-\theta_0\|>\delta_T}m_T(\theta)'\Sigma(\theta,\theta)^{-1}m_T(\theta) \approx CT^{1-2\alpha-2\gamma}$ so if $0<\gamma<1/2-\alpha$, then part (i) of Assumption \ref{ass: strong-id 1} holds. Taylor expansion shows that
$$
\sup_{|\theta-\theta_0|<\delta_T}|m_T(\theta)-M_T(\theta-\theta_0)|\leq T^{1/2-\alpha}q^2\sup_{\theta\in\Theta}\sup_{i,j}\left| \frac{\partial^2m(\theta)}{\partial\theta_i\partial\theta_j}\right|\delta_T^2,
$$
so for $\gamma>1/2(1/2-\alpha)$ part (ii) holds. Finally, $M_T=T^{1/2-\alpha}M$, thus, part (iii) holds if $\gamma<1/2-\alpha$. To summarize, parts (i)-(iii) hold for any $\gamma$ with $1/2(1/2-\alpha)<\gamma<1/2-\alpha$.

\begin{theorem}\label{theorem: strong identification}
Suppose Assumptions  \ref{Ass: Uniform CLT}-\ref{Ass: consistency of covariance} and \ref{ass: strong-id 1} hold, then the $QLR$ statistic defined in equation (\ref{eq: defn of qlr}) converges in distribution to a $\chi^2_q$ uniformly over $\mathcal{P}_0$ as the sample size increases to infinity, while at the same time the conditional critical value  $c_\alpha(h_T,\widehat\Sigma)$ converges in probability to the $1-\alpha$-quantile of a $\chi^2_q$-distribution.
Thus under strong identification the conditional QLR test is asymptotically equivalent to the classical unconditional QLR test under the null.
\end{theorem}

Theorem \ref{theorem: strong identification} concerns behavior under the null but can be extended to local alternatives. Define local alternatives to be sequences of alternatives which are contiguous in the sense of Le Cam (see, for example, chapter 10 in Van der Vaart and Wellner (1996)) with sequences in $\mathcal{P}_0$ satisfying Assumption \ref{ass: strong-id 1}.  By the definition of contiguity, under all such sequences of local alternatives $c_\alpha(h_T,\widehat\Sigma)$ will again converge to a $\chi^2_q$ critical value, implying that our conditional $QLR$ test coincides with the usual $QLR$ test under these sequences.

\section{Concentrating out nuisance parameters}\label{section- concentrating out}

As highlighted in Section \ref{section-setting}, processes $g_T(\cdot)$ satisfying Assumptions \ref{Ass: Uniform CLT}-\ref{Ass: consistency of covariance} arise naturally when considering normalized moment conditions in GMM estimation.  Such processes arise in other contexts as well, however.  In particular, one can often obtain such moment functions by ``concentrating out'' well-identified structural nuisance parameters.  This is of particular interest for empirical work, since in many empirical settings we are interested in testing a hypothesis concerning a subset of the structural parameters, while the remaining structural (nuisance) parameters  are unrestricted.  In this section we show that if we have a well-behaved estimate of the structural nuisance parameters (in a sense made precise below), a normalized moment function based on plugging in this estimator provides a process $g_T(\cdot)$ satisfying Assumptions \ref{Ass: Uniform CLT}-\ref{Ass: consistency of covariance}.  We then show that these results may be applied to test hypotheses on the coefficients on the endogenous regressors in quantile IV models, treating the parameters on the exogenous controls as strongly-identified nuisance parameters.

In this section we assume that we begin with a $(q+p)$-dimensional structural parameter which can be written as $(\beta,\theta)$, where we are interested in testing a hypothesis $H_0:\theta=\theta_0$ concerning only the $q$-dimensional parameter $\theta$.  The hypothesis of interest is thus that there exists some value $\beta_0$ of the nuisance $\beta$ such that the $k$-dimensional moment condition $Eg^{(L)}(\beta_0,\theta_0)=0$  holds.   Here we use superscript $(L)$ to denote the ``long'' or non-concentrated moment condition and define a corresponding ``long'' mean function $m_T ^{(L)}(\beta,\theta)$.  We assume there exists a function $\beta(\theta)$, which we call the pseudo-true value of parameter $\beta$ for a given value of $\theta$, satisfying  $m_T^{(L)}(\beta(\theta_0),\theta_0)=0$.  For values of $\theta$ different from the null value $\theta_0$ the model from which $\beta(\theta)$ comes may be (and often will be) misspecified.  This presents no difficulties for us, as our only requirement will be that there exist an estimator $\widehat\beta(\theta)$ of $\beta(\theta)$ which is $\sqrt{T}$-consistent and asymptotically normal uniformly over $\theta$.  Under additional regularity conditions, we then show that we can use the concentrated moment function $g_T(\theta)=g^{(L)}_T(\theta,\widehat\beta(\theta))$ to implement our inference procedure.

\begin{ass}\label{Ass: concentr 1}
There exists a function $\beta(\theta)$ which for all $\theta$ belongs to the interior of the parameter space for $\beta$  and satisfies $
m_T^{(L)}(\beta(\theta_0),\theta_0)=0,
$ and an estimator $\widehat\beta(\theta)$ such that $\left(g_T^{(L)}(\beta,\theta)-m_T^{(L)}(\beta,\theta), \sqrt{T}(\widehat\beta(\theta)-\beta(\theta))                                                                                \right)$ are jointly uniformly asymptotically normal,
\begin{align*}
\lim_{T\to\infty}\sup_{P\in\mathcal{P}_0}\sup_{f\in BL_1}\left| E_P\left[f\left(
\begin{array}{c}
g_T^{(L)}(\beta,\theta)-m_T^{(L)}(\beta,\theta) \\                                                                                  \sqrt{T}(\widehat\beta(\theta)-\beta(\theta))                                                                                 \end{array}
\right)\right]-E\left[f(\mathbb{G})\right]\right|=0.
\end{align*}
where $\mathbb{G}=(G^{(L)}(\beta,\theta),G_\beta(\theta))$ is a mean-zero Gaussian process with covariance function $\Sigma_L(\beta,\theta,\beta_1,\theta_1)$, such that process $\mathbb{G}$ is uniformly equicontinuous and uniformly bounded over $\mathcal{P}_0$.

\end{ass}
\begin{ass}\label{Ass: concentr 2}
Assume that the covariance function is uniformly bounded, uniformly positive definite, and uniformly continuous in $\beta$ along $\beta(\theta)$.  In particular, for fixed $\bar\lambda>0$ and any sequence $\delta_T\to0$ we have
$$
1/\bar\lambda\leq \inf_{P\in\mathcal{P}_0}\inf_{\theta}\lambda_{\min}(\Sigma_L(\beta(\theta),\theta,\beta(\theta),\theta))\leq \sup_{P\in\mathcal{P}_0}\sup_{\theta}\lambda_{\max}(\Sigma_L(\beta(\theta),\theta,\beta(\theta),\theta))\leq \bar\lambda;
$$
$$
\lim_{T\to\infty}\sup_{P\in\mathcal{P}_0}\sup_{\theta,\theta_1}\sup_{\|\beta-\beta(\theta)\|<\delta_T} \sup_{\|\beta_1-\beta(\theta_1)\|<\delta_T} \|\Sigma_L(\beta,\theta,\beta_1,\theta_1)-\Sigma_L(\beta(\theta),\theta,\beta(\theta_1),\theta_1)\|=0.
$$
\end{ass}

\begin{ass}\label{Ass: concentr 3} There is an estimator $\widehat\Sigma_L(\beta,\theta,\beta_1,\theta_1)$ of $\Sigma_L(\beta,\theta,\beta_1,\theta_1)$ such that
\begin{align*}
\lim_{T\to\infty}\sup_{P\in\mathcal{P}_0}P\left\{\sup_{\beta,\theta,\beta_1,\theta_1}\left\| \widehat\Sigma_L(\beta,\theta,\beta_1,\theta_1)-\Sigma_L(\beta,\theta,\beta_1,\theta_1)\right\| >\varepsilon\right\}=0.
\end{align*}
\end{ass}

\begin{ass}\label{Ass: concentr 4}
For some sequence $\delta_T\to\infty$, $\delta_T/\sqrt{T}\to 0$, for each $P\in\mathcal{P}_0$ there exists a deterministic sequence of $k\times p$ functions $M_T(\theta)$ such that:
\begin{align*}
\lim_{T\to\infty}\sup_{P\in\mathcal{P}_0}\sup_{\theta}\sup_{\sqrt{T}|\beta-\beta(\theta)|\leq \delta_T}\left\| m^{(L)}_T(\beta,\theta)-m^{(L)}_T(\beta(\theta),\theta)-M_T(\theta)\sqrt{T}(\beta-\beta(\theta))\right\|=0.
\end{align*}
We assume that these functions $M_T(\theta)$ are uniformly bounded:
$
\sup_{P\in\mathcal{P}_0}\sup_{\theta}\|M_T(\theta)\|<\infty,
$
and there exists an estimator $\widehat{M}_T(\theta)$ such that
$$
\lim_{T\to\infty}\sup_{P\in\mathcal{P}_0}P\left\{\sup_{\theta}\left\| \widehat{M}_T(\theta)-M_T(\theta)\right\| >\varepsilon\right\}=0.
$$
\end{ass}

\paragraph{Discussion of Assumptions}   Assumptions \ref{Ass: concentr 1}-\ref{Ass: concentr 3} extend Assumptions \ref{Ass: Uniform CLT}- \ref{Ass: consistency of covariance}, adding strong-identification conditions for $\beta$.  In particular, Assumption  \ref{Ass: concentr 1} states that there exists a consistent and asymptotically normal estimator $\widehat\beta(\theta)$ uniformly over $\theta$.   Assumption \ref{Ass: concentr 2} additionally guarantees that the rate of convergence for $\widehat\beta(\theta)$ is uniformly $\sqrt{T}$, and Assumption \ref{Ass: concentr 3} guarantees that the covariance function is well-estimable.  Note that if the estimator  $\widehat\beta(\theta)$ is obtained using some subset of the moment conditions $g^{(L)}$,
the covariance matrix $\Sigma_L$ may be degenerate along some directions, violating Assumption \ref{Ass: concentr 3}.  In such cases we should  reformulate the initial moment condition $g^{(L)}$ to exclude the redundant directions.
Assumption \ref{Ass: concentr 4} supposes that  $m_T$ is linearizable in $\beta$ in the neighborhood of $\beta(\theta)$. In many GMM models  $m_T^{(L)}(\beta,\theta)=\sqrt{T}E\varphi^{(L)}(\beta,\theta)$ and thus we have $$M_T(\theta)=\frac{\partial}{\partial\beta}E\varphi^{(L)}(\beta,\theta)\mid_{\beta=\beta(\theta)}.$$
This last expression is typically consistently estimable provided $E\varphi^{(L)}(X_t,\beta,\theta)$ is twice-continuously-differentiable in $\beta$, in which case Assumption \ref{Ass: concentr 4} comes from Taylor expansion  in $\beta$ around $\beta(\theta)$.  Note the close relationship between Assumption \ref{Ass: concentr 4} and Assumption \ref{ass: strong-id 1} part (ii).

\begin{theorem}\label{thm: concentrating out} Let Assumptions \ref{Ass: concentr 1}-\ref{Ass: concentr 4} hold, then the moment function $g_T(\theta)=g_T^{(L)}(\widehat\beta(\theta),\theta)$, mean function $m_T(\theta)=m_T^{(L)}(\beta(\theta),\theta)$, covariance function
\begin{align*}
\Sigma(\theta,\theta_1)=\left(I_k, M_T(\theta)\right)\Sigma_L(\beta(\theta),\theta,\beta(\theta_1),\theta_1)\left(I_k, M_T(\theta_1)\right)^\prime,
\end{align*}
and its estimate
\begin{align*}
\widehat\Sigma(\theta,\theta_1)=\left(I_k, \widehat{M}_T(\theta)\right)\widehat\Sigma_P(\widehat\beta(\theta), \theta, \widehat\beta(\theta_1),\theta_1)\left(I_k, \widehat{M}_T(\theta_1)\right)^\prime,
\end{align*}
satisfy the Assumptions \ref{Ass: Uniform CLT}-\ref{Ass: consistency of covariance}.
\end{theorem}
The proof of Theorem \ref{thm: concentrating out} may be found in the Supplementary Appendix.

The assumption that the nuisance parameter $\beta$ is strong-identified, specifically the existence of a uniformly-consistent and asymptotically-normal estimator
$\widehat{\beta}(\theta)$ and the linearizability of $m_T^{(L)}(\beta,\theta)$ in $\beta$, plays a key role here. Andrews and Cheng (2012) and Andrews and Mikusheva (2014) show in models with weakly identified nuisance parameters the asymptotic distributions of many statistics will depend on the unknown values of the nuisance parameter, greatly complicating inference.  In such cases, rather than concentrating out the nuisance parameter we may instead use the projection method.  The projection method tests the continuum of hypothesis $H_0:\theta=\theta_0,\beta=\beta_0$ for different values of $\beta_0$, and rejects the null $H_0:\theta=\theta_0$ only if all hypotheses of the form $H_0:\theta=\theta_0,\beta=\beta_0$ are rejected.  Thus, even in cases where the nuisance parameter may be poorly identified one can test $H_0:\theta=\theta_0$ by applying our conditioning method to test a continuum of hypotheses $H_0:\theta=\theta_0,\beta=\beta_0$ provided the corresponding $g_T^{(L)}(\beta,\theta)$ processes satisfy Assumptions \ref{Ass: Uniform CLT}-\ref{Ass: consistency of covariance}.

\subsection{Example: quantile IV regression}\label{section- quantile IV}

To illustrate our results on concentrating out nuisance parameters we consider inference  on the coefficients on the endogenous regressors in a quantile IV model. This setting has been studied in Chernozhukov and Hansen (2008), where the authors used an Anderson-Rubin-type statistic, and in Jun (2008) where K and J statistics were suggested. Here we propose inference based on a QLR statistic.

Consider an instrumental-variables model of quantile treatment effects
as in Chernozhukov and Hansen (2005). Let the data consist of i.i.d. observations on
  an outcome variable $Y_{t}$,  a vector of endogenous regressors $D_{t}$,
 a vector of exogenous controls $C_{t}$, and  a $k\times1$
vector of instruments $Z_{t}$.  Following Chernozhukov and Hansen (2006) we
assume a linear-in-parameters model for the $\tau$-quantile treatment effect,
known up to parameter $\psi=\left(\beta,\theta\right),$ and will
base inference on the moment condition
\begin{align}\label{eq: quantileIV model}
E\left[\left(\tau-\mathbb{I}\left\{ Y_{t}\le C_{t}'\beta_0+D_{t}'\theta_0\right\} \right)\left(\begin{array}{c}
C_{t}\\
Z_{t}
\end{array}\right)\right]=0.
\end{align}

If we were interested in joint inference on the parameters $\left(\beta,\theta\right)$
we could simply view this model as a special case of GMM. In practice, however, we are often concerned with the coefficient
$\theta$ on the endogenous regressor, so  $\beta$ is a nuisance
parameter and we would prefer to conduct inference on $\theta$ alone.
To do this we can follow Jun (2008) and obtain for each value $\theta$
an estimate $\widehat{\beta}\left(\theta\right)$ for $\beta$ by running a standard, linear-quantile regression of $Y_{t}-D_{t}'\theta$
on $C_{t}$. In particular, define
$$
\widehat\beta(\theta)=\arg\min_{\beta}\frac{1}{T}\sum_{t=1}^T \rho_\tau(Y_t- D_t^\prime\theta-C_t^\prime \beta),
$$
where $\rho_\tau(\cdot)$ is the $\tau$-quantile check function.
 The idea of estimating $\widehat\beta(\theta)$ from simple quantile regression, introduced in Chernozhukov and Hansen (2008), is easy to implement and computationally feasible.  Under mild regularity conditions, $\widehat\beta(\theta)$ will be a consistent and asymptotically-normal estimator for the pseudo-true value $\beta(\theta)$ defined by
\begin{align}\label{eq: quantile betatheta}
E\left[\left(\tau-\mathbb{I}\left\{ Y_{t}\le C_{t}'\beta\left(\theta\right)+D_{t}'\theta\right\} \right)C_{t}\right]=0
\end{align}
for each $\theta$. If we then define the concentrated moment function
$$
g_{T}\left(\theta\right)=\frac{1}{\sqrt{T}}\sum_{t=1}^{T}\left(\tau-\mathbb{I}\left\{ Y_{t}\le C_{t}'\widehat{\beta}\left(\theta\right)+D_{t}'\theta\right\} \right)Z_{t},
$$
mean function
$$
m_{T}\left(\theta\right)=\sqrt{T}E\left[\left(\tau-\mathbb{I}\left\{ Y_{t}\le C_{t}'\beta\left(\theta\right)+D_{t}'\theta\right\} \right)Z_{t}\right],
$$
and the covariance estimator
\begin{align*}
\widehat\Sigma(\theta_1,\theta_2)=\frac{1}{T}\sum_{t=1}^T\left[ \left(\tau-\mathbb{I}\{\varepsilon_t(\widehat\beta(\theta_1),\theta_1)<0\}\right)\left( \tau-\mathbb{I}\{\varepsilon_t(\widehat\beta(\theta_2),\theta_2)<0\}\right)\cdot\right.\\ \left.\cdot\left(Z_t-\widehat{A}(\theta_1)C_t\right)\left(Z_t-\widehat{A}(\theta_2)C_t\right)'
\right],
\end{align*}
where $\varepsilon(\beta,\theta)=Y_t- D_t^\prime\theta-C_t^\prime \beta$, $\widehat{A}(\theta)=\widehat{M}_T(\theta)\widehat{J}^{-1}(\theta)$,
\begin{align*}
\widehat{M}_T(\theta) = \frac{1}{Th_T} \sum_{t=1}^T Z_tC_t^\prime k \left(\frac{\varepsilon_t(\widehat\beta(\theta),\theta)}{h_T}\right),
\widehat{J}(\theta) = \frac{1}{Th_T} \sum_{t=1}^T C_tC_t^\prime k \left(\frac{\varepsilon_t(\widehat\beta(\theta),\theta)}{h_T}\right),
\end{align*}
we show in the Supplementary Appendix that these choices satisfy Assumptions \ref{Ass: concentr 1}-\ref{Ass: concentr 4} under the following regularity conditions:

\begin{ass}\label{Ass: quantile IV}
\begin{itemize}
\item[(i)] $(Y_t,C_t,D_t,Z_t)$ are i.i.d., $E\|C\|^{4}+E\|D\|^{2+\varepsilon}+E\|Z\|^{4}$ is uniformly bounded above, and the matrix $E[(C_t',Z_t')(C_t',Z_t')']$ is full rank.
\item[(ii)] The conditional density $f_{\varepsilon(\theta)}(s|C,D,Z)$ of $\varepsilon(\theta)=Y- D^\prime\theta-C^\prime\beta(\theta)$ is uniformly bounded over the support of $(C,D,Z)$ and is twice continuously differentiable at $s=0$  with a second derivative that is  uniformly continuous in $\theta$;
\item[(iii)] For each $\theta$ the value of $\beta(\theta)$ defined in equation (\ref{eq: quantile betatheta}) is in the interior of the parameter space;
\item[(iv)] $\inf_{\theta}\lambda_{\min}(J(\theta))>0$  for $J(\theta)=E\left[f_{\varepsilon(\theta)}(0)CC^\prime\right]$;
\item[(v)]  The kernel $k(v)$ is such that $\sup|k(v)|<\infty,$ $\int|k(v)|dv<\infty,$ $\int k(v)dv=1, $ and $\int k^2(v)dv<\infty$.
\end{itemize}
\end{ass}

Under Assumption \ref{Ass: quantile IV}, one may use the $QLR$ statistic paired with conditional critical values to construct confidence sets for $\theta$ in this model.  In Section \ref{section- numerical} we provide simulation results comparing the performance of QLR tests with known alternatives.  Both Chernozhukov and Hansen (2008) and Jun (2008) suggested  Anderson-Rubin type statistics for this model which have stable power, but which are inefficient in overidentified models under strong identification. To overcome this inefficiency, Jun (2008) introduced a K test analogous to Kleibergen (2005).  This test is locally efficient under strong identification and has good power for small violations of the null hypothesis regardless of identification strength.  However, K tests often suffer from substantial declines in power at distant alternatives.  To overcome this deficiency a number of approaches to combining the K and AR statistics have been suggested by different authors, including the JK test discussed by Jun (2008), which is expected to improve power against distant alternatives but is inefficient under strong identification.  By contrast, our approach allows one to use QLR tests, which retain efficiency under strong identification without sacrificing power at distant alternatives.

\section{Numerical performance of the conditional QLR test}\label{section- numerical}
In this section we examine the performance of the conditional QLR
test in two numerical examples, first simulating the performance of
the conditional QLR test in a quantile IV model and then constructing
confidence sets for Euler equation parameters in US data by inverting
the conditional QLR test.

\subsection{Simulations: quantile IV  model}

We simulate the performance of the QLR test in a quantile IV model with a single endogenous regressor and $k$ instruments. We draw i.i.d. random vectors $\left(U_{t},D_{t},Z_{t}'\right)'=(\Phi^{-1}(\xi_{U,t}),\Phi^{-1}(\xi_{D,t}),\Phi^{-1}(\xi_{Z_1,t}),...,\Phi^{-1}(\xi_{Z_k,t}))$ from
a Gaussian copula. In particular, the $\xi$'s are normals with mean zero, all  variances equal to  one, $\mbox{cov}(\xi_U,\xi_D)=\rho$, $\mbox{cov}(\xi_D,\xi_{Z_j})=\pi$ and all other covariances are zero, and $\Phi$ is the standard-normal distribution function.
We generate the outcome variable $Y_{t}$ from the location-scale
model,
\[
Y_{t}=\gamma_{1}+\gamma_{2}D_{t}+\left(\gamma_{3}+\gamma_{4}D_{t}\right)\left(U_{t}-\frac{1}{2}\right),
\]
which  implies a linear conditional-quantile model for all quantiles. The only control variable, $C_t$, is a constant.
For our simulations we focus on the median, $\tau=\frac{1}{2}$ and the corresponding coefficients are $\beta=\gamma_1$ and $\theta=\gamma_2$.

In this model, we can think of $\rho$ as measuring the endogeneity
of the regressor $D_{t}$: if $\rho=0$ then there is no endogeneity
and a linear quantile regression of $Y_{t}$ on $D_{t}$ and a constant
will yield consistent estimates of $(\beta,\theta)$. If
on the other hand $\rho\neq0$, we need to adopt a quantile IV strategy to obtain consistent estimates.
The parameter ${\pi}$ controls the strength of the identification
under the quantile IV approach,
so the model will be partially identified when ${\pi}=0$ and
weakly identified when ${\pi}$ is close to zero.

We are interested in inference on the coefficient $\theta$ on the endogenous
regressor, treating the intercept $\beta$ as a nuisance
parameter and calculating our conditional QLR test as described in
Section \ref{section- quantile IV}. For comparison we also calculate the weak-instrument-robust AR, K, and JK tests of Jun (2008), which are based
on the same concentrated moment conditions but use different test statistics.  In Jun (2008)'s simulations the test suggested by Chernozhukov and Hansen (2008) performed quite similarly to Jun's AR test, so here we report results only for Jun's tests.

\subsubsection{Simulation results}

Our simulations set $\gamma_{i}=1$ for all $i$ so the true value
of our coefficient of interest $\theta$ is 1. We fix $\rho=0.25$
and consider samples of 1,000 observations generated from the model
above as we vary the identification parameter ${\pi}$. We
considered cases with five and ten instruments, $k=5$ and $k=10$,
but for brevity here report only the results for five instruments:
the results for ten instruments are quite similar and are available
upon request.

Table \ref{tab:QIV test size 5 instruments} reports the simulated
size of nominal 5\% tests for the null $H_{0}:\theta=1$ as we vary
the identification parameter ${\pi}$. As we would hope given
the identification-robust nature of the tests studied, the simulated
size is in all cases close to the nominal level 5\% and is insensitive
to the strength of identification as measured by ${\pi}$.

\begin{table}
\begin{tabular}{|c|c|c|c|c|c|c|c|c|}
\hline
 $\pi $& 0.02 & 0.04 & 0.06 & 0.08&0.1& 0.15 & 0.25 & 0.4\tabularnewline
\hline
\hline
AR & 5.09\% & 5.25\%& 5.15\% & 5.04\% & 5.09\% & 5.00\% & 5.26\% & 5.18\%\tabularnewline
\hline
K & 5.64\% & 5.16\% & 5.14\% & 5.13\% & 5.46\% & 4.98\% & 4.87\% & 5.17\%\tabularnewline
\hline
JK & 5.27\% & 5.25\% & 5.39\% & 5.05\% & 5.43\% & 5.14\% & 5.05\% & 5.46\%\tabularnewline
\hline
QLR & 5.62\% & 5.12\% & 5.18\% & 5.06\% & 4.99\% & 5.04\% & 5.22\% & 5.18\%\tabularnewline
\hline
\end{tabular}

\caption{Power nominal 5\% tests in quantile IV simulations with five instruments
and 1,000 observations. Based on 10,000 simulation replications, and
10,000 draws of conditional critical values. \label{tab:QIV test size 5 instruments}}
\end{table}

\begin{figure}
  \begin{sideways}
    \begin{minipage}{23.5cm}
     \includegraphics[scale=0.51]{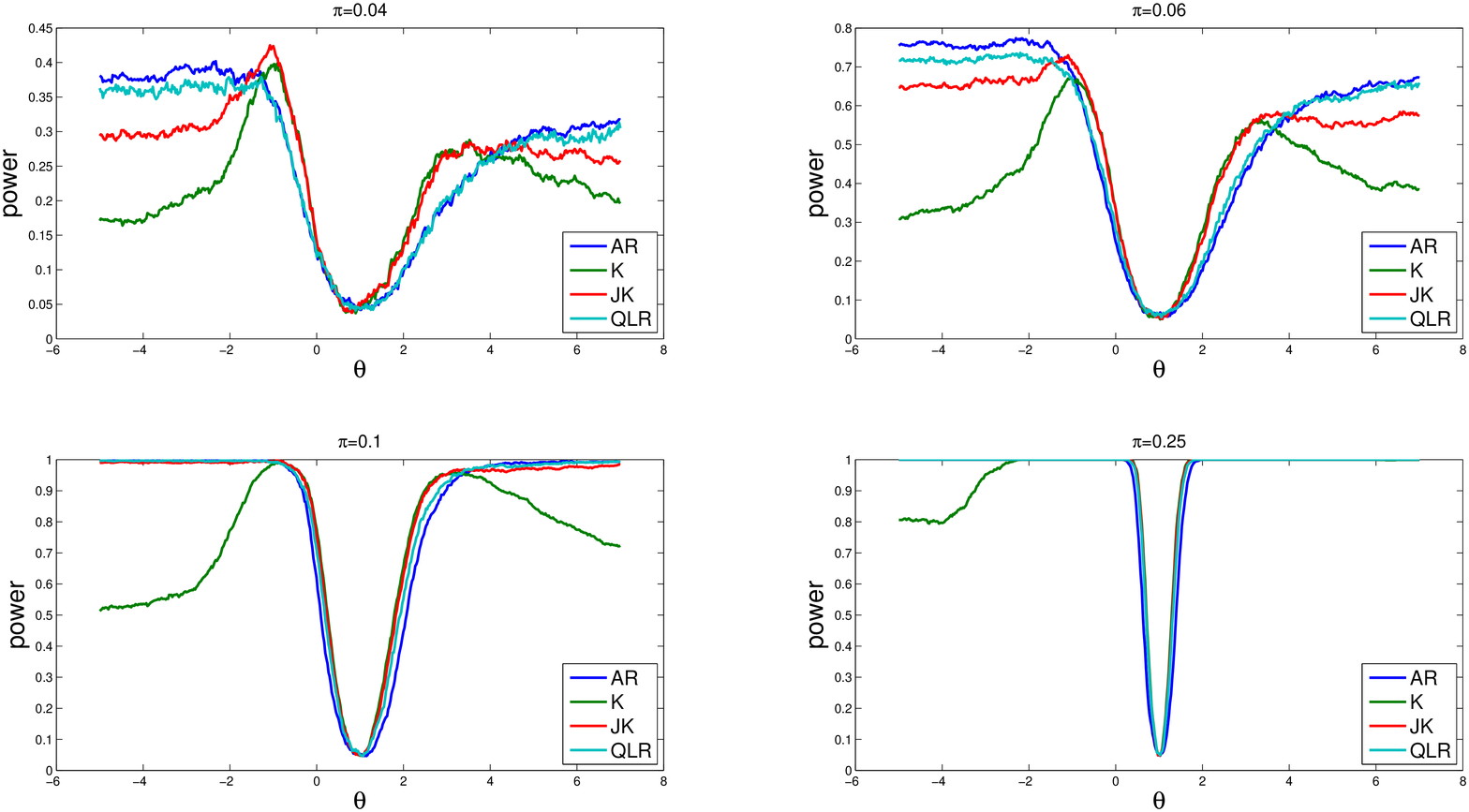}
    \end{minipage}
  \end{sideways}
  \centering
  \caption{Power of nominal 5\% tests in quantile IV simulations with five instruments,
1,000 observations, and four different values of identification strength
${\pi}$. Based on 1,000 simulation replications and 10,000 draws
of conditional critical values. \label{fig:QIV test power 5 instruments-1}}
\end{figure}

Since all tests considered have approximately correct size, we next
compare them in terms of power. Figure \ref{fig:QIV test power 5 instruments-1}
plots the simulated power
of the tests for a range of values for the identification
strength parameter ${\pi}$. Since the scale of Figure \ref{fig:QIV test power 5 instruments-1} makes the power curves difficult
to distinguish in the well-identified cases, Figure \ref{fig:QIV test power 5 instruments, Zoomed rho_pi=00003D0.4}
plots power curves for ${\pi}=0.4$ focusing on a smaller
neighborhood of the null.

\begin{figure}
\includegraphics[scale=0.9]{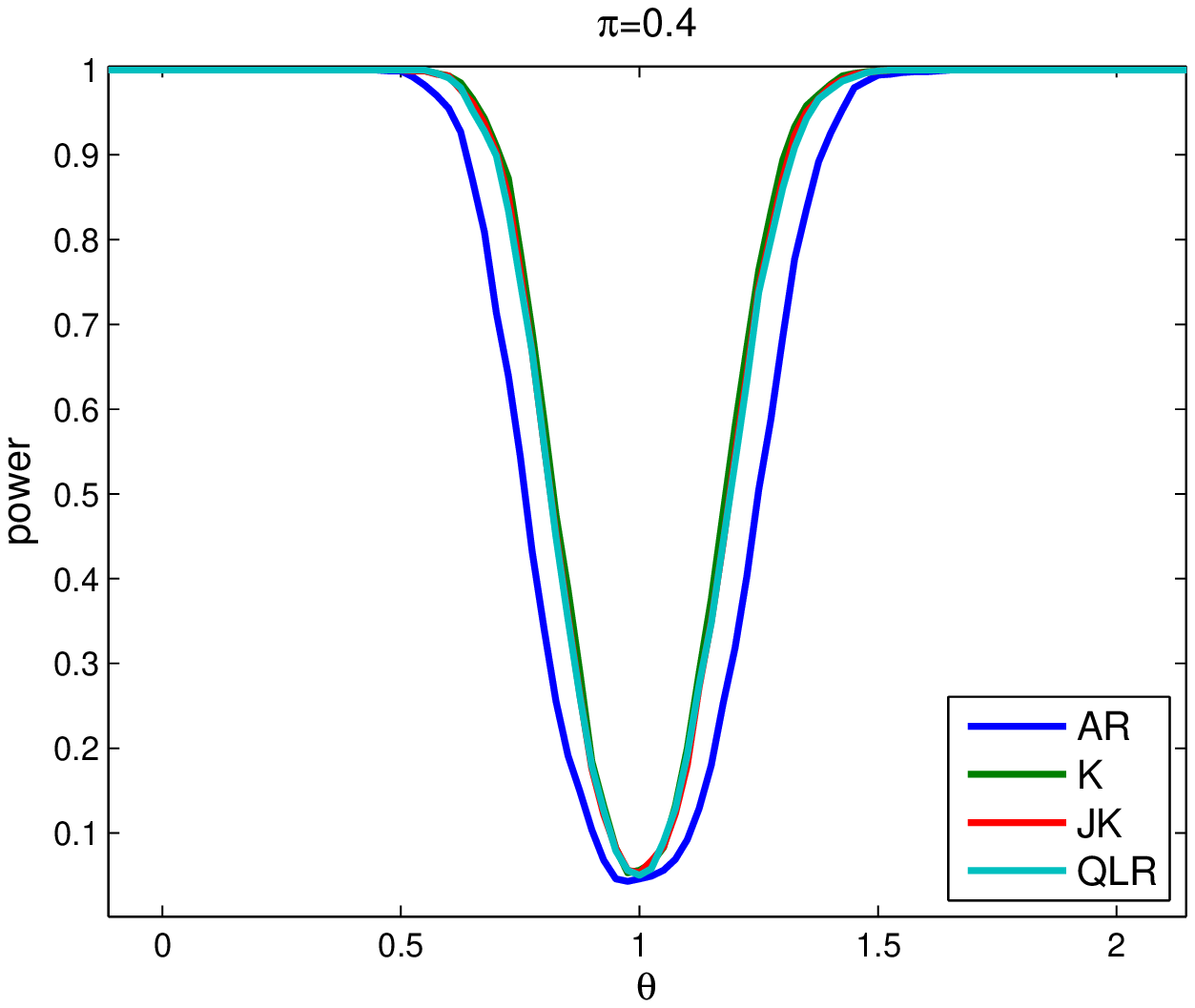}

\caption{Power nominal 5\% tests in quantile IV simulations with five instruments,
1,000 observations, and ${\pi}=0.4$. Based on 1,000 simulation
replications and 10,000 draws of conditional critical values. \label{fig:QIV test power 5 instruments, Zoomed rho_pi=00003D0.4}}
\end{figure}

From these figures we can see that when identification is quite weak
(that is, when ${\pi}$ is close to zero), all tests have
power substantially below one. The K and JK tests tend to have good
power close to the null but often suffer from substantial declines
in power as one moves away from the null. By contrast, the power of
the AR and QLR tests generally tends to increase as we consider alternatives
more distant from the null.
For ${\pi}$ large the power curves of the K, JK, and QLR tests
are essentially indistinguishable local to the null $\theta=1$, while
the AR test is clearly inefficient in this case. Despite its good power
close to the null we see that even in this case the K test continues
to exhibit pronounced power deficiencies against some alternatives,
consistent with the results of Jun (2008). If we fix ${\pi}\neq0$
and take the sample size to infinity the K and (by the results of
Theorem \ref{theorem: strong identification}) QLR tests will be efficient local to $\theta=1$.
By contrast, the JK and AR tests will be inefficient, though the degree
of inefficiency for the JK test will be small. Thus, we see that the
conditional QLR test we propose has appealing power properties; it is efficient when identification is strong and does not experience power declines at distant alternatives when identification is weak.

\subsection{Empirical example: Euler equation}
\label{subsection- Euler Equation}

As an empirical example, we invert the QLR and several other robust
tests to calculate identification-robust confidence sets based on
the nonlinear Euler equation specification discussed in Example 1.
Following Stock and Wright (2000) we use an extension of the long
annual data-set of Campbell and Shiller (1987). Our specification
corresponds to the CRRA-1 specification of Stock and Wright (2000),
which takes $C_{t}$ to be aggregate consumption, $R_{t}$ to be an
aggregate stock market return and $Z_{t}$ to contain of a constant,
$C_{t-1}/C_{t-2},$ and $R_{t-1}$, resulting in a three-dimensional
moment condition ($k=3$)- see Stock and Wright (2000) for details.
As in Kleibergen (2005), to estimate all covariance matrices we use
the Newey-West estimator with one lag.%
\footnote{While the model implies that $\sqrt{T}g_{T}(\cdot)$ is a martingale when evaluated
at the true parameter value, the QLR statistic also depends on the
behavior of $g_{T}$ away from the null. Likewise, Kleibergen (2005)
notes the importance of using a HAC covariance matrix estimator in
the construction of the K statistic. We could use a martingale-difference
covariance estimator in constructing the $S$ statistic, but doing
so substantially increases the size of the joint $S$ confidence set
for $\left(\delta,\gamma\right)$ so we focus on the HAC formulation
for comparability with the other confidence sets studied.%
} We first construct a confidence set for the full parameter vector
$\theta=\left(\delta,\gamma\right)$ and then consider inference on
the risk-aversion coefficient $\gamma$ alone.

\subsubsection{Confidence sets for the full parameter vector}

\begin{figure}
  \begin{sideways}
    \begin{minipage}{23.5cm}
     \includegraphics[scale=0.51]{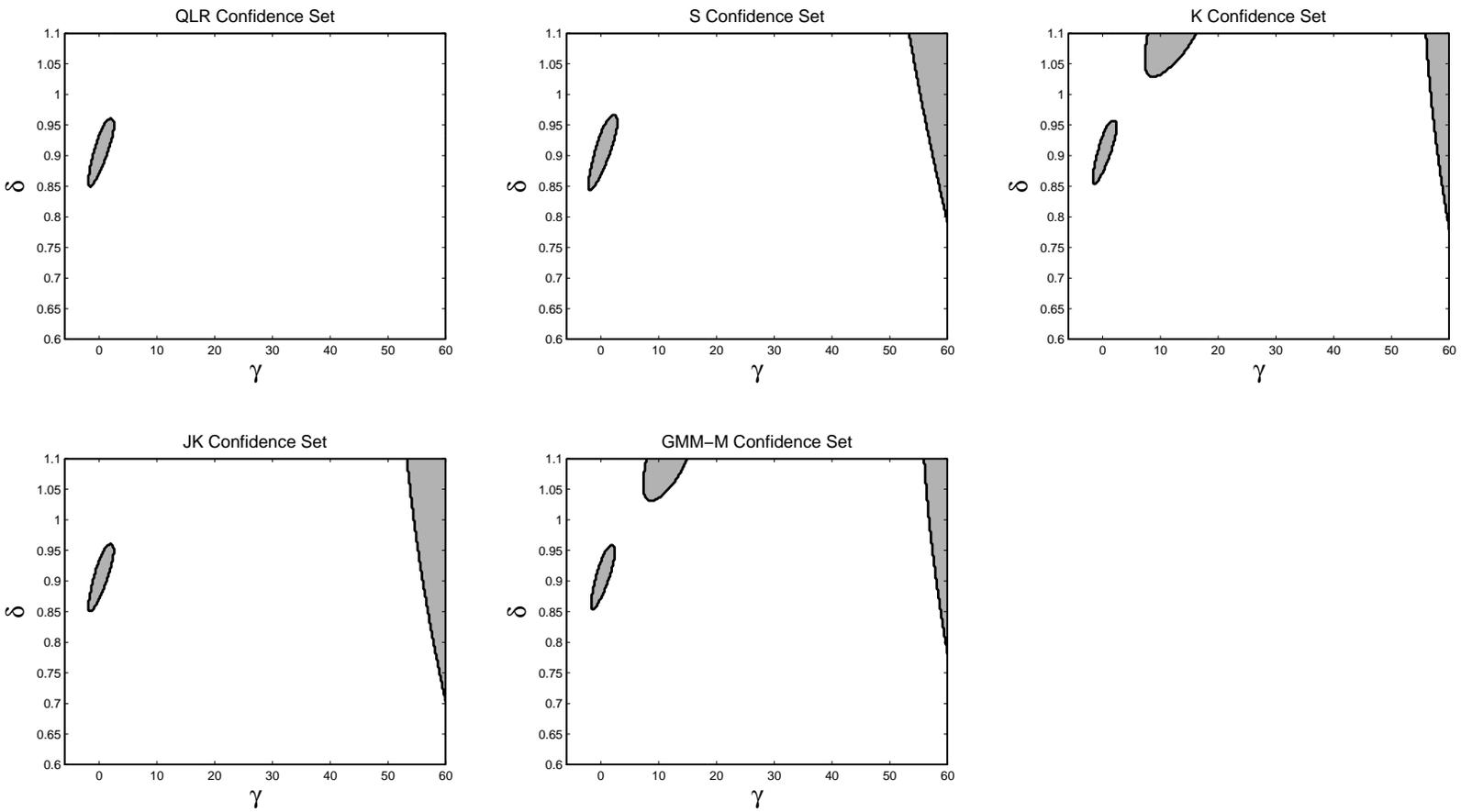}
    \end{minipage}
  \end{sideways}
  \centering
  \caption{Joint 90\% QLR, S, K, JK, and QCLR confidence sets for risk aversion
($\gamma)$ and the discount factor $(\delta)$ based on annual data,
three moment conditions, and 1,000 draws of critical values.\label{fig:Joint non-QLR CS}}
\end{figure}

Joint 90\% confidence sets for $\theta=\left(\delta,\gamma\right)$
based on inverting QLR, S, K, JK, and GMM-M tests of Stock and Wright (2000)
and Kleibergen (2005) are reported in Figure \ref{fig:Joint non-QLR CS}.%
\footnote{Note that our S confidence set differs from that of Stock and Wright
(2000) which, in addition to assuming that the summands in $g_{T}(\theta_0)$ are serially uncorrelated, also assumes conditional homoskedasticity.%
} As we can see, the QLR confidence set is substantially smaller than
the others considered, largely by virtue of eliminating disconnected
components of the confidence set. To quantify this difference, note
that the S, K, JK, and QCLR confidence sets cover 4.3\%, 4.43\%, 5.46\%,
and 4.5\% of the parameter space $\left(\delta,\gamma\right)\in\left[0.6,1.1\right]\times\left[-6,60\right]$,
respectively, while the QLR confidence set covers only 0.64\% of the
parameter space.

\subsection{Confidence sets for risk aversion}

Stock and Wright (2000) argued that once one fixes the risk-aversion
parameter $\gamma$ the discount factor $\delta$ is well identified.
Under this assumption we calculate conditional QLR confidence sets
for $\gamma$ based on two approaches, first by plugging in an estimator
for $\delta$ based on the moment condition instrumented with a constant
and then concentrating out $\delta$ using the continuous-updating
estimator (CUE), where in each case we modify the moment conditions
as discussed in Section \ref{section- concentrating out} to account for this estimation.
For comparison we consider the S, K, JK, and GMM-M tests evaluated
at the restricted CUE for $\delta$ which, as Stock and Wright (2000)
and Kleibergen (2005) argue, allow valid inference under the assumption
that $\delta$ is well identified. The resulting confidence sets are
reported in Table \ref{tab:Gamma confidence sets}. Unlike in the
joint confidence set case we see that the QLR confidence set is larger
than the JK confidence set, but it is nonetheless the second smallest
confidence set out of the five considered. Further, we see that in
this application concentrating out the nuisance parameter using the
CUE results in a smaller confidence set than does plugging in the
estimate based on the moment condition instrumented with a constant.

\begin{table}
\begin{tabular}{|c|c|c|}
\hline
 & 90\% Confidence Set & Length\tabularnewline
\hline
\hline
QLR- Constant Instrument & $\left[-2,1.7\right]$ & 3.7\tabularnewline
\hline
QLR- CUE & $\left[-1.3,1.9\right]$ & 3.2\tabularnewline
\hline
S & $\left[-1.6,2.3\right]$ & 3.9\tabularnewline
\hline
K & $\left[-1.1,1.8\right]\cup\left[8,12.3\right]$ & 7.2\tabularnewline
\hline
JK & $\left[-1.2,1.9\right]$ & 3.1\tabularnewline
\hline
GMM-M & $\left[-1.1,1.8\right]\cup\left[8,12.3\right]$ & 7.2\tabularnewline
\hline
\end{tabular}

\caption{90\% confidence sets for risk aversion parameter $\gamma$, treating
nuisance parameter $\delta$ as well identified, based on annual data.\label{tab:Gamma confidence sets}}
\end{table}

\section{Conclusions}\label{section: conclude}

This paper argues that moment-equality models without any identification assumptions have a functional nuisance parameter.  We introduce a sufficient statistic for this nuisance parameter and construct conditional tests.  Our results substantially expand the set of statistics available in weakly- or partially-identified models, and in particular allow the use of quasi-likelihood ratio statistics, which often have  superior power properties compared to the widely-used Anderson-Rubin type statistics.  We show that our tests have uniformly correct asymptotic size over a large class of models, and find that the proposed tests perform well in simulations in a quantile IV model and give smaller confidence sets than existing alternatives in a nonlinear Euler equation model.

\section{References}

\noindent\hspace{-10mm}   Andrews, D.W.K. and X. Cheng (2012):
``Estimation and Inference with Weak, Semi-strong and Strong
Identification,'' \emph{Econometrica}, 80(5), 2153-2211.

\noindent\hspace{-10mm}   Andrews, D.W.K., X. Cheng, and P. Guggenberger (2011):
``Generic Results for Establishing the Asymptotic Size of Confidence Sets and Tests,'' \emph{unpublished manuscript}.

\noindent\hspace{-10mm}  Andrews D.W.K. and P. Guggenberger (2009): ``Hybrid and Size-Corrected Subsampling Methods,'' \emph{Econometrica}, 77, 721-762.

\noindent\hspace{-10mm}  Andrews D.W.K. and P. Guggenberger (2010): ``Asymptotic Size and a Problem with Subsampling and with the $m$ out of $n$ Bootstrap,'' \emph{Econometric Theory}, 26, 426-468.

\noindent\hspace{-10mm} Andrews, D.W.K., M. Moreira, and J. Stock
(2006): ``Optimal Two-Sided Invariant Similar Tests for Instrumental
Variables Regression,'' \emph{Econometrica}, 74, 715-752.

\noindent\hspace{-10mm} Andrews, D.W.K., M. Moreira, and J. Stock
(2008): ``Efficient Two-sided Nonsimilar Invariant Tests in IV Regression with Weak Instruments,'' \emph{Journal of Econometrics}, 146, 241-254.

\noindent\hspace{-10mm}  Andrews, I. and A. Mikusheva (2014):
``A Geometric Approach to Weakly Identified Econometric Models,'' \emph{unpublished manuscript}

\noindent\hspace{-10mm}  Campbell, J.Y. and R.J. Shiller (1987): ``Cointegration Tests of
Present Value Models,'' \emph{Journal of Political Economy}, 95,
1062-1088.

\noindent\hspace{-10mm}  Chernozhukov, V. and C. Hansen (2005): ``An IV Model of Quantile
Treatment Effects,'' \emph{Econometrica}, 73, 245-261.

\noindent\hspace{-10mm}  Chernozhukov, V. and C. Hansen (2006): ``Instrumental Quantile Regression
Inference for Structural and Treatment Effect Models,'' \emph{Journal
of Econometrics}, 132, 491-525.

\noindent\hspace{-10mm}  Chernozhukov, V. and C. Hansen (2008): ``Instrumental Variable Quantile Regression: A Robust Inference Approach,'' \emph{Journal of Econometrics}, 142, 379-398.

\noindent\hspace{-10mm}  Dedecker, J. and S. Louhichi (2002): ``Maximal Inequalities and Empirical Central Limit Theorems,'' in  H. Dehling, T. Mikosch and M. Sorensen (eds.) \emph{Empirical Process Techniques for Dependent Data}, 137-161.

\noindent\hspace{-10mm}  Hansen, L.P. (1982) : ``Large Sample Properties of Generalized Method of Moments Estimators,'' \emph{Econometrica}, 50, 1029-1054.

\noindent\hspace{-10mm}  Hansen, L.P. and K. Singleton (1982):``Generalized Instrumental Variables Estimation of Nonlinear Rational Expectations Models,'' \emph{Econometrica}, 50, 1269-1286.

\noindent\hspace{-10mm}  Jun, S.J. (2008): ``Weak Identification Robust Tests in an Instrumental
Quantile Model,'' \emph{Journal of Econometrics}, 144, 118-138.

\noindent\hspace{-10mm}  Kleibergen, F. (2005): ``Testing Parameters in GMM without Assuming
that They are Identified,'' \emph{Econometrica}, 73, 1103-1124.

\noindent\hspace{-10mm}  Lehmann, E.L. and J.P. Romano (2005): \emph{Testing Statistical Hypotheses}, New York: Springer; 3rd edition.

\noindent\hspace{-10mm}  McFadden, D. (1989): ``A Method of Simulated Moments for Estimation
of Discrete Response Models without Numerical Integration,'' \emph{Econometrica},
57, 995-1026.

\noindent\hspace{-10mm} Moreira, M. (2003): ``A Conditional
Likelihood Ratio Test for Structural Models,'' \emph{Econometrica},
71, 1027-1048.

\noindent\hspace{-10mm}  Pakes, A. and D. Pollard (1989): ``Simulation and the Asymptotics
of Optimization Estimators,'' \emph{Econometrica}, 57, 1027-1057.

\noindent\hspace{-10mm}  Schennach, S. (2014): ``Entropic Latent Variable Integration by Simulation,''
\emph{Econometrica}, 82, 345-385.

\noindent\hspace{-10mm}  Staiger, D. and J. Stock (1997): ``Instrumental Variables Regression
with Weak Instruments,'' \emph{Econometrica}, 65, 557-586.

\noindent\hspace{-10mm}  Stock, J. and J. Wright (2000): ``GMM with Weak Identification,''
\emph{Econometrica}, 82, 345-385.

\noindent\hspace{-10mm}  Van der Vaart, A.W. and J.A. Wellner (1996): \emph{Weak Convergence and
Empirical Processes}. New York: Springer.

\section{Appendix}

\paragraph{Proof of Lemma \ref{lemma- sufficient statistic}.}
The proof trivially follows from equation (\ref{eq: g as a function of h}) and the observations that (i) the distribution of $g_T(\theta_0)\sim N(0,\Sigma(\theta_0,\theta_0))$ does not depend on $m_T(\cdot)$, (ii) the function $\Sigma(\theta,\theta_0)\Sigma(\theta_0,\theta_0)^{-1}$ is deterministic and known, and (iii) the vector $g_T(\theta_0)$ is independent of $h_T(\cdot)$. $\Box$

\paragraph{Proof of Theorem \ref{Theorem: uniform asymptotics}.}
Let us introduce the process
$$
G_h(\theta)=H(G,\Sigma)(\theta)=G(\theta)-\Sigma\left(\theta,\theta_{0}\right) \Sigma\left(\theta_{0},\theta_{0}\right)^{-1}G\left(\theta_{0}\right),
$$
and a random variable $\xi=G(\theta_0)$ which is independent of $G_h(\cdot)$.
First, we notice that Assumptions \ref{Ass: Uniform CLT}-\ref{Ass: consistency of covariance} imply that
$\eta_T=(g_T(\theta_0),h_T(\cdot)-m_T(\cdot),\widehat\Sigma(\cdot,\cdot))$ converges uniformly to $\eta=(\xi,G_h(\cdot),\Sigma(\cdot,\cdot))$, that is,
\begin{align}\label{eq: in the proof of asymp}
\lim_{T\to\infty}\sup_{P\in\mathcal{P}_{0}} \sup_{f\in BL_{1}} \left|E_{P}\left[f\left(\eta_T\right)\right] -E\left[f\left(\eta\right)\right]\right|=0,
\end{align}
where $BL_1$ is again the class of bounded Lipshitz functionals with  constant 1.
We assume here that the distance on the space of realizations is measured as follows: for $\eta_i=(\xi_i,G_{h,i}(\cdot),\Sigma_i(\cdot,\cdot))$ (for $i=1,2$),  $$d(\eta_1,\eta_2)=\|\xi_1-\xi_2\|+\sup_\theta\|G_{h,1}(\theta)-G_{h,2}(\theta)\|+ \sup_{\theta,\tilde\theta}\|\Sigma_1(\theta,\tilde\theta)-\Sigma_2(\theta,\tilde\theta)\|.$$
Statement (\ref{eq: in the proof of asymp}) then follows from the observation that the function which takes $(G(\cdot),\Sigma(\cdot,\cdot))$ to $(\xi,G_h(\cdot), \Sigma(\cdot,\cdot))$ is Lipshitz in $(G,\Sigma)$ if $|\xi|<C$ for some constant $C$, provided $\Sigma$ satisfies Assumption \ref{Ass: Bounded Covariance}.

Next, note that  for $\varsigma_T=(g_T(\theta_0),h_T(\cdot),\widehat\Sigma(\cdot,\cdot))$ and $\tilde\varsigma_T=(\xi,G_h(\cdot)+m_T,\Sigma(\cdot,\cdot))$ we have
\begin{align}\label{eq: in the proof of asymp1}
\lim_{T\to\infty}\sup_{P\in\mathcal{P}_{0}} \sup_{f\in BL_{1}} \left|E_{P}\left[f\left(\varsigma_T\right)\right] -E\left[f\left(\tilde\varsigma_T\right)\right]\right|=0,
\end{align}
as follows from (\ref{eq: in the proof of asymp}) and the observation that bounded Lipshitz functionals of $\varsigma_T$ are also bounded Lipshitz in $\eta_T$.

Let us introduce the function $F(x)=\mathbb{I}\{x<C_1\}+\frac{C_2-x}{C_2-C_1}\mathbb{I}\{C_1\leq x<C_2\}$ for some $0<C_1<C_2$ and consider the functional
$$
R_C(\xi, h ,\Sigma)=R(\xi,h,\Sigma)F(\xi'\Sigma\left(\theta_{0},\theta_{0}\right)^{-1}\xi),
$$
which is a continuous truncation of the functional $R(\xi,h,\Sigma)=R(g,\Sigma)$. Consider the conditional quantile function corresponding to the new statistic
$$
c_{C,\alpha}(h,\Sigma) =\inf \left\{ c: P^*\left\{\xi:R_C(\xi,h,\Sigma)\leq c\right\}\geq 1-\alpha \right\}.
$$
As our next step we show that
$c_{C,\alpha}(h,\Sigma)$ is Lipshitz in $h(\cdot)$ and $\Sigma(\cdot,\cdot)$ for all $h$ with $h(\theta_0)=0$ and $\Sigma$ satisfying Assumption \ref{Ass: Bounded Covariance}.

Assumption \ref{Ass: statistic} implies that there exists a constant $K$ such that
$$
\left\|R_C(\xi,h_1,\Sigma)-R_C(\xi,h_2,\Sigma)\right\|\leq Kd(h_1,h_2)
$$
for all $\xi, h_1, h_2$ and $\Sigma$. Let $c_i=c_{C,\alpha}(h_i,\Sigma)$, then
$$
1-\alpha\leq P^*\left\{\xi:R_C(\xi,h_1,\Sigma)\leq c_1\right\}\leq P^*\left\{\xi:R_C(\xi,h_2,\Sigma)\leq c_1+Kd(h_1,h_2)\right\}.
$$
Thus $c_2\leq c_1+Kd(h_1,h_2)$. Analogously we get $c_1\leq c_2+Kd(h_1,h_2)$, implying that $c_{C,\alpha}$ is Lipshitz in $h$.
The same argument shows that $c_{C,\alpha}$ is Lipshitz in $\Sigma$.

Assume the conclusion of Theorem \ref{Theorem: uniform asymptotics} does not hold. Then there exists some $\delta>0$, an infinitely increasing sequence of sample sizes $T_i$, and a sequence of probability measures $P_{T_{i}}\in\mathcal{P}_0$ such that for all $i$
$$
P_{T_{i}}\left\{ R(g_{T_{i}}(\theta_0),h_{{T_{i}}},\widehat\Sigma)> c_\alpha (h_{T_{i}},\widehat\Sigma)+\varepsilon\right\}>\alpha+\delta.
$$
Choose $C_1$ such that
$$
\limsup P_{T_{i}}\left\{ g_{T_{i}}(\theta_0)'\widehat\Sigma\left(\theta_{0},\theta_{0}\right)^{-1}g_{T_{i}}(\theta_0)\geq C_1 \right\}<\frac{\delta}{2},$$
which can always be done since according to Assumption \ref{Ass: Uniform CLT} $g_T(\theta_0)$ converges uniformly to $N(0,\Sigma\left(\theta_{0},\theta_{0}\right))$. Since
\begin{align*}
P_T \left\{ R> x\right\} \leq P_T\left\{ R_C> x \right\}+P_T\left\{ g_{T}(\theta_0)'\widehat\Sigma\left(\theta_{0},\theta_{0}\right)^{-1}g_{T}(\theta_0)\geq C_1\right\},
\end{align*}
and $c_{C,\alpha} (h_T,\widehat\Sigma)< c_\alpha (h_T,\widehat\Sigma)$ we have that for all $i$
\begin{align}\label{eq: eq2}
P_{T_{i}}\left\{ R_C(g_{T_{i}}(\theta_0),h_{{T_{i}}},\widehat\Sigma)\geq c_{C,\alpha} (h_{{T_{i}}},\widehat\Sigma)+\varepsilon\right\}>\alpha+\frac{\delta}{2}.
\end{align}
Denote by $\mathcal{T}_T$ a random variable distributed as   $R_C(\xi_T,h_T,\widehat\Sigma)- c_{C,\alpha} (h_T,\widehat\Sigma)$ under the law $P_T$, and by $\mathcal{T}_{\infty,T}$ a random variable distributed as   $R_C(\xi,G_{h}+m_T,\Sigma)- c_{C,\alpha} (G_{h}+m_T,\Sigma)$ under the law $P_T$.  The difference between these variables is that the first uses the finite-sample distribution of $(\xi_T, h_T,\widehat\Sigma)$, while the latter uses their asymptotic counterparts  $(\xi,G_h+m_T,\Sigma)$. Equation (\ref{eq: in the proof of asymp1}) and the bounded Lipshitz property of the statistic $R_C$ and the conditional critical value imply that
\begin{align}\label{eq:  eq1}
\lim_{T\to\infty}\sup_{f\in BL_1}\left|Ef(\mathcal{T}_T)-Ef(\mathcal{T}_{\infty,T})\right|=0.
\end{align}
Since $\mathcal{T}_{T_{i}}$ is a sequence of bounded random variables, by Prokhorov's theorem  there exists a subsequence $T_j$ and a random variable $\mathcal{T}$ such that $\mathcal{T}_{T_j}\Rightarrow\mathcal{T}$. By (\ref{eq:  eq1}), $\mathcal{T}_{\infty,T_j}\Rightarrow\mathcal{T}$. Since (\ref{eq: eq2}) can be written as $P\{\mathcal{T}_T\geq \varepsilon\}>\alpha+\delta/2$,
$$
\lim\inf P\{\mathcal{T}_{\infty,T_j}>0\}\geq P\{\mathcal{T}>0\}\geq P\{\mathcal{T}\geq\varepsilon\}\geq\lim\sup P\{\mathcal{T}_{T_j}\geq\varepsilon\}\ge\alpha+\frac{\delta}{2}.
$$
However, from the definition of quantiles we have
$$
P\{\mathcal{T}_{\infty,T_j}>0\}= P_T\left\{ R_C(\xi,G_{h}+m_T,\Sigma)> c_{C,\alpha} (G_{h}+m_T,\Sigma)\right\}\leq \alpha,
$$
since the statistic $\mathcal{T}_{\infty,T_j}$ is the statistic in the limit problem and so controls size by Lemma \ref{lemma- sufficient statistic}. Thus we have reached a contradiction. $\Box$

\paragraph{Proof of Theorem \ref{theorem: strong identification}.}
As shown in Theorem \ref{Theorem: uniform asymptotics}, Assumptions \ref{Ass: Uniform CLT}-\ref{Ass: consistency of covariance} imply that the distribution of the QLR statistic is uniformly asymptotically approximated by the distribution of the same statistic in the limit problem. Thus, it suffices to prove the statement of Theorem \ref{theorem: strong identification} for the limit problem only, which is to say when $g_T(\cdot)$ is Gaussian process with mean $m_T(\cdot)$ and known covariance $\Sigma$. In our case $QLR=R(g_T(\theta_0),h_T,\Sigma)$, where
\begin{align}
R(\xi,h,\Sigma)=\xi'\Sigma(\theta_0,\theta_0)^{-1}\xi-\inf_{\theta}\left(V(\theta)\xi+h(\theta)\right)' \Sigma(\theta,\theta)^{-1}\left(V(\theta)\xi+h(\theta)\right),\label{eq: note 3}
\end{align}
and $V(\theta)= \Sigma\left(\theta,\theta_{0}\right) \Sigma\left(\theta_{0},\theta_{0}\right)^{-1}$.
Denote by $\mathcal{A}$ the event $\mathcal{A}=\{g_T(\theta_0)'\Sigma(\theta_0,\theta_0)^{-1}g_T(\theta_0)<C\}$ and note that by choosing the constant $C>0$ large enough we can guarantee that the probability of $\mathcal{A}$ is arbitrarily close to one.

Let $\hat\theta_T$ be  the value at which the optimum in (\ref{eq: note 3}) is achieved (the case when the optimum may not be achieved may be handled similarly, albeit with additional notation). We first show that $\hat\theta_T\to^p\theta_0$. For any $a,b$ we have  $(a+b)^2\geq \frac{a^2}{2}-b^2$, so
\begin{align}\label{eq: note 1}
(V(\theta)g_T(\theta_0)+h_T(\theta))'\Sigma(\theta,\theta)^{-1}(V(\theta)g_T(\theta_0)+h_T(\theta))\geq \frac{1}{2} m_T(\theta)'\Sigma(\theta,\theta)^{-1}m_T(\theta)\\- (V(\theta)g_T(\theta_0)+h_T(\theta)-m_T(\theta) )'\Sigma(\theta,\theta)^{-1}(V(\theta)g_T(\theta_0)+h_T(\theta)-m_T(\theta)).\notag
\end{align}
Assumptions \ref{Ass: Bounded Covariance} and \ref{ass: strong-id 1} (vi) guarantee that the second term on the right-hand side of (\ref{eq: note 1}) is stochastically bounded, so denote this term $A(\theta)$.  For any probability $\varepsilon>0$ there exists a constant $C$ such that
\begin{align*}
\inf_{P\in\mathcal{P}_0}P\left\{\sup_{\theta\in\Theta}A(\theta)\leq C \mbox{ and } \mathcal{A}\right\}\geq 1-\varepsilon.
\end{align*}
Assumption \ref{ass: strong-id 1}(i) implies that there exists $T_1$ such that for all $T>T_1$ and $P\in\mathcal{P}_0$ we have
$$
\inf_{ \|\theta-\theta_0\|>\delta_T}m_T(\theta)'\Sigma(\theta,\theta)^{-1}m_T(\theta)>4C.
$$
Putting the last three inequalities together we get that for $T>T_1$ and all $P\in\mathcal{P}_0$
\begin{align*}
P\left\{\inf_{ \|\theta-\theta_0\|>\delta_T}(V(\theta)g_T(\theta_0)+h_T(\theta))' \Sigma(\theta,\theta)^{-1}(V(\theta)g_T(\theta_0)+h_T(\theta))>C \mbox{ and } \mathcal{A}\right\}\geq 1-\varepsilon.
\end{align*}
This implies that $\sup_{P\in\mathcal{P}_0}P\left\{\|\hat\theta_T-\theta_0\|>\delta_T\right\}\leq\varepsilon$  for all  $T>T_1.$

As our second step we show that for any $\varepsilon>0$
\begin{align}\label{eq: note 4}
\lim_{T\to\infty}\sup_{P\in\mathcal{P}_0}P\left\{\left|\inf_{ \|\theta-\theta_0\|<\delta_T}g_T(\theta)' \Sigma(\theta,\theta)^{-1}g_T(\theta)\right.\right.\\
\left.\left.-\inf_{ \|\theta-\theta_0\|<\delta_T}\tilde{g}_T (\theta)'\Sigma(\theta_0,\theta_0)^{-1}\tilde{g}_T(\theta)\right|>\varepsilon\right\}=0\notag,
\end{align}
where we replace the process $g_T(\theta)=V(\theta)g_T(\theta_0)+h_T(\theta)$ by the process $\tilde{g}_T(\theta)=g_T(\theta_0)+m_T(\theta)$ with the same mean function $m_T(\theta)$ and covariance $\tilde{\Sigma}(\theta,\theta_1)=\Sigma(\theta_0,\theta_0)$ for all $\theta,\theta_1$. For this new process we have $\tilde{V}(\theta)=I$ and $\tilde{h}_T(\theta)=m_T(\theta)$.
To verify (\ref{eq: note 4}), restrict attention to the event $\mathcal{A}$ for some large $C>0$. The functional that transforms $(g_T(\theta_0),h,\Sigma(\theta,\theta),V(\cdot))$ to $\inf_{ \|\theta-\theta_0\|<\delta_T}(V(\theta)g_T(\theta_0)+h(\theta))'\Sigma(\theta,\theta)^{-1}( V(\theta)g_T(\theta_0)+h(\theta))$
is Lipshitz  in $h$, $V$ and $\Sigma(\theta,\theta)$ on $\mathcal{A}$. Thus,
\begin{align*}
\left|\inf_{ \|\theta-\theta_0\|<\delta_T}g_T(\theta)' \Sigma(\theta,\theta)^{-1}g_T(\theta)-\inf_{ \|\theta-\theta_0\|<\delta_T}\tilde{g}_T (\theta)'\Sigma(\theta_0,\theta_0)^{-1}\tilde{g}_T(\theta)\right|\\ \leq K_1 \sup_{\|\theta-\theta_0\|\leq \delta_T}|h_T(\theta)-m_T(\theta)|+K_2 \sup_{\|\theta-\theta_0\|\leq \delta_T}\left\|\Sigma(\theta,\theta)-\Sigma(\theta_0,\theta_0)\right\|\\+K_3\sup_{\|\theta-\theta_0\|\leq \delta_T}\left\|\Sigma(\theta,\theta_0)-\Sigma(\theta_0,\theta_0)\right\|.
\end{align*}
Note, however, that $h_T(\theta)-m_T(\theta)=G(\theta)-\Sigma\left(\theta,\theta_{0}\right) \Sigma\left(\theta_{0},\theta_{0}\right)^{-1}G(\theta_0)$. Assumptions \ref{ass: strong-id 1} (iv) and (v) therefore imply (\ref{eq: note 4}).

As our third step, we linearly approximate $m_T$ using Assumption \ref{ass: strong-id 1} (ii), which implies that for any $\varepsilon>0$
\begin{align*}
\lim_{T\to\infty}\sup_{P\in\mathcal{P}_0}P\left\{\left|\inf_{ \|\theta-\theta_0\|<\delta_T}(g_T(\theta_0) +m_T(\theta))'\Sigma(\theta_0,\theta_0)^{-1}(g_T(\theta_0)+m_T(\theta))\right.\right.\\ \left.\left.
-\inf_{ \|\theta-\theta_0\|<\delta_T}(g_T(\theta_0) +M_T(\theta-\theta_0))'\Sigma(\theta_0,\theta_0)^{-1}(g_T(\theta_0)+M_T(\theta-\theta_0))\right|>\varepsilon\right\}=0.
\end{align*}
Indeed, on the set $\mathcal{A}$ we have that $\inf_{ \|\theta-\theta_0\|<\delta_T}(g_T(\theta_0) +m(\theta))'\Sigma(\theta_0,\theta_0)^{-1}(g_T(\theta_0)+m(\theta))$ is Lipshitz in $m$.

So far we have shown that
$QLR$ is asymptotically equivalent to
$$
QLR_1=g_T(\theta_0)'\Sigma(\theta_0,\theta_0)^{-1}g_T(\theta_0)- \inf_{ \|\theta-\theta_0\|<\delta_T}(g_T(\theta_0) +M_T(\theta-\theta_0))'\Sigma(\theta_0,\theta_0)^{-1}(g_T(\theta_0)+M_T(\theta-\theta_0)),
$$
and in particular that $QLR-QLR_1\to^p0$ as $T\to \infty$.
Note, however, that statistic
$$
QLR_2=g_T(\theta_0)'\Sigma(\theta_0,\theta_0)^{-1}g_T(\theta_0)- \inf_{ \theta}(g_T(\theta_0) +M_T(\theta-\theta_0))'\Sigma(\theta_0,\theta_0)^{-1}(g_T(\theta_0)+M_T(\theta-\theta_0))
$$
is $\chi^2_{q}$ distributed provided $M_T$ is full rank. The difference between $QLR_1$ and $QLR_2$ is in the area of optimization, and the optimizer in $QLR_2$ is $$
\theta^*=(M_T'\Sigma(\theta_0,\theta_0)^{-1}M_T)^{-1}M_T'\Sigma(\theta_0,\theta_0)^{-1}g_T(\theta_0)\sim N(0,(M_T'\Sigma(\theta_0,\theta_0)^{-1}M_T)^{-1}).$$
Assumption \ref{ass: strong-id 1} (iii) guarantees that $\|\theta^*\|/\delta_T$ converges uniformly to zero in probability, and thus that
$$
\lim_{T\to\infty}\sup_{P\in\mathcal{P}_0}P\{\|\theta^*-\theta_0\|>\delta_T\}=0.
$$
As a result,
$QLR_1-QLR_2\to^p0$, which proves that $QLR\Rightarrow\chi^2_{q}$ uniformly over $\mathcal{P}_0.$
The convergence of the conditional critical values is proved in a similar way.$\Box$

\end{document}